%% file: ms.tex
\documentclass[12pt]{amsart}
\usepackage{amsfonts, amsmath, amssymb, amsthm, bibtex, caption, commath, dsfont, enumitem, graphicx, lscape, mathdots, mathrsfs, mathtools, mdframed, multicol, pst-plot, tikz, ulem, xcolor}

\usepackage[makeroom]{cancel}
\usepackage[headings]{fullpage}

\usetikzlibrary{external}
\usepackage{hyperref}
\hypersetup{
    colorlinks=true,
    linkcolor=blue,    
    urlcolor=blue,
    citecolor=blue
}

\everymath{\displaystyle} % for math
\setenumerate[1]{label=(\arabic*)}

\newcommand{\N}{\mathds{N}}
\newcommand{\Q}{\mathds{Q}}

\newcommand{\Z}{\mathds{Z}}

\newcommand{\oo}{\infty}
\newcommand{\epts}{:}

\renewcommand{\phi}{\varphi}
\renewcommand{\epsilon}{\varepsilon}

\renewcommand{\emph}{\textsl}

\DeclareMathOperator{\spn}{span}

\DeclareMathOperator{\acyc}{acyc}
\DeclareMathOperator{\sink}{sink}
\DeclareMathOperator{\st}{\vert}
\DeclareMathOperator{\sgn}{sgn}
\DeclareMathOperator{\proj}{Proj}

\newcommand{\ncr}[2]{\mathchoice
{\genfrac () {0pt} 0 {#1} {#2}}
{\genfrac () {0pt} 1 {#1} {#2}}
{\genfrac () {0pt} 2 {#1} {#2}}
{\genfrac () {0pt} 3 {#1} {#2}}}

\newtheorem{theorem}{Theorem}[section]
\newtheorem{lemma}[theorem]{Lemma}
\newtheorem*{theorem*}{Theorem}

\newtheorem{corollary}[theorem]{Corollary}
\theoremstyle{definition}
\newtheorem{example}[theorem]{Example}
\newtheorem{definition}[theorem]{Definition}

\newtheorem*{observation}{Observation}

\tikzstyle{vertex}=[fill, circle, radius = 2pt]
\tikzstyle{small}=[fill, circle, inner sep = 2pt]
\tikzstyle{colored}=[draw=black, circle, minimum size = 0.8cm]

\usepackage{xparse} % for evil screwed up commands

\NewDocumentCommand{\pstv}{m m O{-} O{above} O{sloped}}{\draw [line width=0.5mm, #3] (#1) -- (#2) node [midway, #4, #5] {$+$}}

\NewDocumentCommand{\ltr}{m m O{-} O{above} O{sloped}}{\draw [line width=0.5mm, #3] (#1) -- (#2) node [midway, #4, #5] {$+$} node [near start, sloped] {\tikz \draw[-latex, sloped, line width=0.75mm] (0.3,0) -- (0.6,0);} node [near end, sloped] {\tikz \draw[-latex, sloped, line width=0.75mm] (0.3,0) -- (0.6,0);}}

\NewDocumentCommand{\rtl}{m m O{-} O{above} O{sloped}}{\draw [line width=0.5mm, #3] (#1) -- (#2) node [midway, #4, #5] {$+$} node [near start, sloped] {\tikz \draw[-latex, sloped, line width=0.75mm] (0.6,0) -- (0.3,0);} node [near end, sloped] {\tikz \draw[-latex, sloped, line width=0.75mm] (0.6,0) -- (0.3,0);}}

\NewDocumentCommand{\ngtv}{m m O{-} O{above} O{sloped}}{\draw [line width=0.5mm, #3] (#1) -- (#2) node [midway, #4, #5] {$-$}}

\NewDocumentCommand{\ptn}{m m O{-} O{above} O{sloped}}{\draw [line width=0.5mm, #3] (#1) -- (#2) node [midway, #4, #5] {$-$} node [near start, sloped] {\tikz \draw[-latex, sloped, line width=0.75mm] (0.3,0) -- (0.6,0);} node [near end, sloped] {\tikz \draw[-latex, sloped, line width=0.75mm] (0.6,0) -- (0.3,0);}}

\NewDocumentCommand{\ntp}{m m O{-} O{above} O{sloped}}{\draw [line width=0.5mm, #3] (#1) -- (#2) node [midway, #4, #5] {$-$} node [near start, sloped] {\tikz \draw[-latex, sloped, line width=0.75mm] (0.6,0) -- (0.3,0);} node [near end, sloped] {\tikz \draw[-latex, sloped, line width=0.75mm] (0.3,0) -- (0.6,0);}}

\NewDocumentCommand{\edge}{m m O{-}}{\draw [line width=0.5mm, #3] (#1) -- (#2)}

\NewDocumentCommand{\drct}{m m m O{-} O{midway}}{\draw [line width=0.5mm, #4] (#1) -- (#2) node [#5, sloped] {\tikz \draw[#3, sloped, line width=0.75mm] (0.3,0) -- (0.6,0);}}

\title{An Extension of Stanley's Symmetric Acyclicity Theorem to Signed Graphs}
\author[O. Coppola]{Oscar Coppola}
\author[J. Huryn]{Jake Huryn}
\author[M. Reilly]{Michael Reilly}

\begin{document}
\maketitle
\begin{abstract}
    In 1995, Richard Stanley introduced the chromatic symmetric function $X_G$ of a graph $G$ and proved that, when written in terms of the elementary symmetric functions, it reveals the number of acyclic orientations of $G$ with a given number of sinks. In this paper, we generalize this result to signed graphs, that is, to graphs whose edges are labeled with $+$ or $-$ and whose colorings and orientations can interact with their signs.
    
    Additionally, we introduce a non-homogeneous basis which detects the number of sinks and which not only gives a Stanley-type result for signed graphs, but gives an analogous result of this form for unsigned graphs as well.
    
    % In 1982, Thomas Zaslavsky proved that the chromatic polynomial $\chi_\Sigma$ of a signed graph $\Sigma$ evaluated at $-1$ equals the number of acyclic orientations of $\Sigma$, generalizing a result of Richard Stanley's. Stanley generalized his own result in 1995 by devising the symmetric chromatic function $X_G$ of a graph $G$, which reveals the number of acyclic orientations of $G$ with any given number of sinks.
    
    % We finally generalize Stanley's 1995 result to signed graphs. We defined the $B$-symmetric chromatic function $X_\Sigma$ of a signed graph $\Sigma$ and prove a theorem analogous to Stanley's. In doing so, we also find a basis for the space of all $B$-symmetric functions, the augmented elementary $B$-symmetric basis.
\end{abstract}

\input{Sections/Introduction.tex}

\input{Sections/Signed_Graphs.tex}

\input{Sections/Main_Theorems.tex}

\input{Sections/deletion_contraction.tex}

\input{Sections/covering.tex}

\input{Sections/linear_extensions.tex}

\input{Sections/sink_counting_appendix.tex}

\input{Sections/phi_calcs.tex}

\input{Sections/main_theorem_proofs.tex}

\input{Sections/further_results.tex}

\input{Sections/bib.tex}

\end{document}

%% file: Sections/Introduction.tex
% \begin{enumerate}
%     \item Chromatic polynomial
%     \item symmetric generalization
%     \item signed generalization
%     \item describe main results
%     \item Acknowledgements 
% \end{enumerate}
\section{Introduction}
A fundamental notion in the study of graphs is that of a \emph{proper coloring} of a graph. This is a function which assigns a natural number to each vertex of a graph in such a way that no two vertices which are connected by an edge are assigned the same color. The \emph{chromatic polynomial} $\chi_G(\lambda)$ of a graph $G$, is a polynomial in $\lambda$ whose value is the number of ways to properly color a graph $G$ in $\lambda$ colors. In 1973, Richard Stanley \cite{stanley73} proved that $\chi_G(-1)$ is the number of acyclic orientations of $G$ (up to sign), a surprising fact which extracts non-coloring information from the chromatic polynomial.

Stanley generalized his result in 1995 \cite{stanley95} by defining the \emph{chromatic symmetric} function $X_G$ of a graph $G$. If $\mathcal P(G)$ is the set of proper colorings of $G$ in colors from $\N$, then
\[
    X_G(x_1,x_2,\dots) = \sum_{\kappa\in\mathcal P(G)}x_{\kappa(v_1)}x_{\kappa(v_2)}\cdots x_{\kappa(v_n)}
\]
where $v_1,v_2,\dots,v_n$ are the of vertices of $G$. Note that this makes $X_G$ a formal power series over countably many variables $\{x_k\st k\in\N\}$. $X_G$ is called symmetric because for any permutation $\pi:\N\to\N$, we may observe that $X_G(x_{\pi(1)},x_{\pi(2)},\dots)=X_G(x_1,x_2,\dots)$. 
%Symmetric functions can be expressed in the \emph{elementary symmetric basis}, consisting of $e_n:= \sum_{1\leq i_1<\cdots<i_n}x_{i_1}\cdots x_{i_n}$ for $n\in \N$
This allows us to state Stanley's result \cite[Theorem 3.3]{stanley95}
\begin{theorem*}
    If the chromatic symmetric function of a graph $G$ is written in terms of the elementary symmetric basis, then the number of acyclic orientations of $G$ with $k$ sinks is the sum of the coefficients of terms having $k$ elementary symmetric function factors.
\end{theorem*}

Generalizing Stanley's work to \emph{signed graphs}, Thomas Zaslavsky introduced the \emph{signed chromatic polynomial} of a signed graph $\Sigma$ and analogously proved that, when evaluated at $-1$, it returns the number of acyclic orientations of $\Sigma$ \cite[Corollary 4.1]{zaslavsky}. In this paper, we will generalize both Zaslavsky's theorem on the signed chromatic polynomial and Stanley's theorem on the chromatic symmetric function to a result about the \emph{chromatic $B$-Symmetric function}, defined as follows.
 
Following E.\ Egge \cite{Egge}, we define the chromatic $B$-symmetric function $X_\Sigma$ of a signed graph $\Sigma$ in direct analogy to the chromatic symmetric function, to be
\[
    X_\Sigma(\dots,x_{-1},x_0,x_1,\dots) = \sum_{\kappa\in\mathcal P(\Sigma)}x_{\kappa(v_1)}x_{\kappa(v_2)}\cdots x_{\kappa(v_n)}
\]
where once again $v_1,v_2,\dots,v_n$ are the of vertices of $\Sigma$, and $\mathcal P(\Sigma)$ is the set of proper colorings of $\Sigma$, which we fully define in the next section. Any permutation $\pi:\Z\to\Z$ satisfying $\pi(k)=-\pi(-k)$ will fix $X_\Sigma$, that is,
\[
    X_\Sigma(\dots,x_{\pi(-1)}, x_{\pi(0)},x_{\pi(1)},\dots)=X_\Sigma(\dots,x_{-1},x_0,x_1,\dots).
\]
Note that, differing from \cite{Egge}, we allow the vertex color $0$, which is in line with Zaslavsky's work in \cite{zaslavsky}. In particular, for any permutation $\pi$  which satisfies $\pi(k)=-\pi(-k)$, we must have $\pi(0)=0$, which distinguishes $0$ from the other colors. This group of permutations is isomorphic to the Coxeter group of type $B$; hence we call it the ``chromatic $B$''-symmetric function, and generally define a function in $\Z$-indexed variables to be $B$-symmetric when it is fixed by all such permutations.

We will define a basis for $B$-symmetric functions which satisfies a Stanley-type result and we will call this the \emph{augmented elementary B-symmetric basis}. Specifically, we have following main theorem.
\begin{theorem*}
    If the chromatic $B$-symmetric function of a signed graph $\Sigma$ is written in terms of the augmented elementary $B$-symmetric basis, then the number of acyclic orientations of $\Sigma$ with $k$ sinks is the sum of the coefficients of terms having $k$ elementary $B$-symmetric function factors.
\end{theorem*}

We will first define all necessary information about signed graphs, explore the properties of chromatic $B$-symmetric functions, and then introduce important tools that we use to prove useful results about basis sets and our main theorem.

\subsection*{Acknowledgements}
This work was done as part of the Summer 2020 and Summer 2021 Knots and Graphs undergraduate research program at The Ohio State University.
The authors were supported financially by the OSU Honors Program Research fund.
We warmly thank Sergei Chmutov, who organized the program and provided generous non-financial support, and Richard Stanley for helpful conversation.
Finally, we are indebted to past years' participants who studied the chromatic ($B$-)symmetric function, including James Enouen, Eric Fawcett, Kat Husar, Hannah Johnson, Rushil Raghavan, and Ishaan Shah.

%% file: Sections/Signed_Graphs.tex
\section{Signed Graphs}
% \begin{enumerate}
%     \item Zaslavsky's stuff (signed graphs, colorings, orientations, covering graphs)
%     \item B-Symmetric Chromatic (also weighted contraction deletion)
%     \item Orientation preserving colorings
%     \item linear extensions
%     \item B-symm = sum over orientations of sum over linear extensions
% \end{enumerate}
As defined by T. Zaslavsky \cite{zaslavsky}, a \emph{signed graph} is a graph whose edges have been labeled with either a plus sign or a minus sign (such edges are called positive or negative respectively). More specifically, we may think of a signed graph as a graph along with a function, $\sgn$, which sends the edges of the graph to an element of the set $\{+,-\}$.

\begin{figure}[h!t]
    \includegraphics{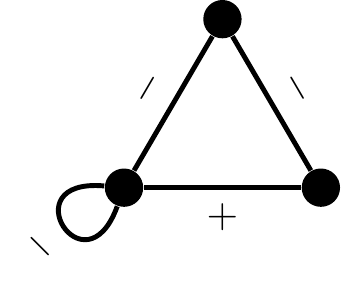}
    \caption{A signed graph with three vertices, three negative edges (one of which is a loop) and one positive edge.}
    % \begin{tikzpicture}
    %     \node[vertex] (v1) at (0,0) {};
    %     \node[vertex] (v2) at (2,0) {};
    %     \node[vertex] (v3) at (1,1.712) {};

    %     \pstv {v1} {v2} [] [below];
    %     \ngtv {v1} {v3} [] [above] [sloped];
    %     \ngtv {v2} {v3} [] [above] [sloped];
    %     \draw[line width=0.5mm, loop below, out = 175, in = 250, looseness = 10] (v1) to (v1);
    %     % \ltr {v1} {v2} [] [below];
    %     % \ntp {v1} {v3} [] [above] [sloped];
    %     % \ptn {v2} {v3} [] [above] [sloped];
    %     % \draw[latex-latex, line width=0.5mm, loop below, out = 175, in = 250, looseness = 10] (v1) to (v1);
    %     \node[rotate = -45] (sign) at (-0.85,-0.6) {$-$};
    % \end{tikzpicture}
\end{figure}

An \emph{orientation} of a signed graph is a way of assigning an arrow to each half-edge (or equivalently, an arrow to each incidence between a vertex and an edge). Each arrow can point either toward or away from its vertex, subject to the condition that on a positive edge both arrows must point in the same direction and on a negative edge both arrows must point in opposite directions.

\begin{figure}[h!t]
    \includegraphics{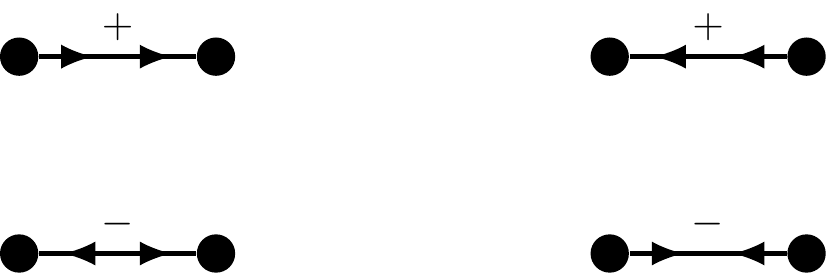}
    \caption{The four possible orientation of an edge.}
    % \begin{tikzpicture}
    % \node[vertex] (+1) at (0,0) {};
    % \node[vertex] (+2) at (2,0) {};
    % \node[vertex] (+3) at (6,0) {};
    % \node[vertex] (+4) at (8,0) {};
    % \node[vertex] (-1) at (0,-2) {};
    % \node[vertex] (-2) at (2,-2) {};
    % \node[vertex] (-3) at (6,-2) {};
    % \node[vertex] (-4) at (8,-2) {};
    
    % \ltr {+1} {+2} [] [above];
    % \rtl {+3} {+4} [] [above];
    % \ntp {-1} {-2} [] [above];
    % \ptn {-3} {-4} [] [above]; 
    % \end{tikzpicture}
\end{figure}

\begin{figure}[h!b]
    \includegraphics{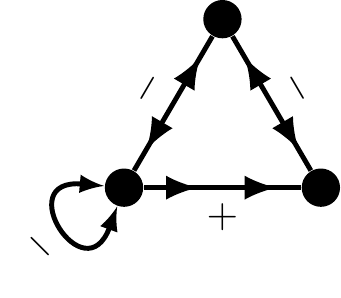}
    \caption{An oriented signed graph.}
    % \begin{tikzpicture}
    %     \node[vertex] (v1) at (0,0) {};
    %     \node[vertex] (v2) at (2,0) {};
    %     \node[vertex] (v3) at (1,1.712) {};

    %     \ltr {v1} {v2} [] [below];
    %     \ntp {v1} {v3} [] [above] [sloped];
    %     \ptn {v2} {v3} [] [above] [sloped];
    %     \draw[latex-latex, line width=0.5mm, loop below, out = 175, in = 250, looseness = 10] (v1) to (v1);
    %     \node[rotate = -45] (sign) at (-0.85,-0.6) {$-$};
    % \end{tikzpicture}
\end{figure}

\begin{figure}[b]\label{cyclic_and_acyclic}
    \includegraphics{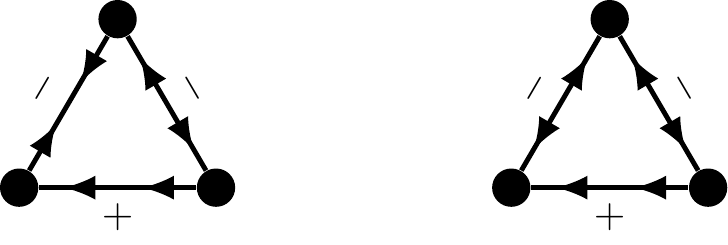}
    \caption{The orientation on the left is cyclic whereas the orientation on the right is acyclic. The acyclic orientation on the right has two sinks: the top vertex and the bottom left vertex.}
    % \begin{tikzpicture}
    %     \node[vertex] (v1) at (0,0) {};
    %     \node[vertex] (v2) at (2,0) {};
    %     \node[vertex] (v3) at (1,1.712) {};

    %     \rtl {v1} {v2} [] [below];
    %     \ptn {v1} {v3} [] [above] [sloped];
    %     \ptn {v2} {v3} [] [above] [sloped];
        
    %     \node[vertex] (u1) at (5,0) {};
    %     \node[vertex] (u2) at (7,0) {};
    %     \node[vertex] (u3) at (6,1.712) {};

    %     \rtl {u1} {u2} [] [below];
    %     \ntp {u1} {u3} [] [above] [sloped];
    %     \ptn {u2} {u3} [] [above] [sloped];
        
    % \end{tikzpicture}
\end{figure}

In this paper, we will be primarily concerned with \emph{acyclic} orientations. A cycle is a closed walk such that every vertex in the walk has at least one arrow pointing into it and one arrow pointing out of it, when only the edges in the walk are considered. See figure 4. An orientation is \emph{acyclic} if it contains no cycles.

When a vertex has only arrows pointing toward it, it is called a \emph{sink}. If it has only arrows pointing away from it, it is called a \emph{source}.

We will also be considering \emph{colorings} of signed graphs. A coloring of a signed graph $\Sigma$ is a function from the vertex set of $\Sigma$ to $\Z$. Using the notation $e\epts uv$ to denote that the endpoints of an edge $e$ are the vertices $u$ and $v$, we say that a coloring $\kappa$ is \emph{proper} if for any $e\epts uv$, we have $\kappa(u) \neq \sgn(e) \cdot \kappa(v)$. Thus if two vertices are connected by a positive edge then a proper coloring cannot assign them the same integer, and if they are connected by a negative edge then a proper coloring cannot assign them integers which are negatives of each other.
%$\kappa(v)$ is usually referred to as the $\emph{color}$ of the vertex $v$.

Some more quick remarks can be made on the nature of signed graph colorings. If a signed graph $\Sigma$ has a vertex with a positive loop, then it has no proper coloring. On the other hand, vertices with negative loops can be properly colored by any integer excluding zero. Finally, we note that a signed graph with all positive edges is virtually indistinguishable from an unsigned graph, and its chromatic $B$-symmetric function is, by re-indexing, identical to its chromatic symmetric function.

%% file: Sections/Main_Theorems.tex
\section{Main Theorems}

Recall that as defined above, the chromatic B-symmetric function of a signed graph $\Sigma$ is
\[
    X_\Sigma(\dots,x_{-1},x_0,x_1,\dots) = \sum_{\kappa\in\mathcal P(\Sigma)}x_{\kappa(v_1)}x_{\kappa(v_2)}\cdots x_{\kappa(v_n)}
\]
where $v_1,v_2,\dots, v_n$ are the vertices of $\Sigma$ and $\mathcal{P}(\Sigma)$ is the set of proper colorings of $\Sigma$. For notational convenience, we will put $x^{\kappa} = x_{\kappa(v_1)}x_{\kappa(v_2)}\cdots x_{\kappa(v_n)}$ whenever applicable.

We will use the notation $p_{a,b} \coloneqq \sum_{i\in \Z}x_{i}^ax_{-i}^b$ and refer to the set $\{p_{a,b}\st a\geq 1, b\geq 0\}\cup \{x_0\}$ as the \emph{p-basis} for the set of B-symmetric functions.
%\section{Main results}
% We know that a\ny $B$-symmetric chromatic function can be written uniquely in terms of sums and products of elements of the set $\{p_{a,b}\st a\geq 1, b\geq 0\}\cup \{x_0\}$. Define $q_{a,b} = (-1)^{a+b+1}p_{a,b}$ and $z = -x_0$. 
It turns out that $X_{\Sigma}$ can always be written uniquely in terms of sums and products of elements of the set $\{p_{a,b}\st a\geq 1, b\geq 0\}\cup \{x_0\}$. We will refer to this set as the \emph{p-basis}.

Defining the \emph{elementary B-symmetric functions} in the natural way,
$$
e_n \coloneqq \sum_{i_1< i_2 < \dots < i_n}x_{i_1}x_{i_2}\cdots x_{i_n}
$$
we can rewrite Newton's identities in the form $p_{n,0} = (-1)^{n+1}ne_n+\sum_{i=1}^{n-1}(-1)^{n+i-1}e_{n-i}p_{i,0}$. 

For convenience we will put $q_{a,b} = (-1)^{a+b+1}p_{a,b}$ and $z = -x_0$. Now we have that any chromatic $B$-symmetric function can be written uniquely in terms of sums an products of elements of the set $\{e_n\st n\geq 1\}\cup \{q_{a,b}\st a,b\geq 1\}\cup \{z\}$. We will call this the \emph{augmented elementary $B$-symmetric basis}. Now we return to our main theorem, whose proof we will postpone until later.
\begin{theorem}\label{main_1}    If the chromatic $B$-symmetric function of a signed graph $\Sigma$ is written in terms of the augmented elementary $B$-symmetric basis, then the number of acyclic orientations of $\Sigma$ with $k$ sinks is the sum of the coefficients of terms having $k$ elementary $B$-symmetric function factors, $e_n$.
\end{theorem}
\begin{example}\label{main_thm_1_ex}
    Let $\Sigma$ be the signed graph
    $\begin{tikzpicture}[baseline = 0.25cm]
        \node[small] (v1) at (0,0) {};
        \node[small] (v2) at (0.75,0) {};
        \node[small] (v3) at (0.375,0.65) {};
        \pstv {v1} {v2} [] [below];
        \ngtv {v1} {v3};
        \pstv {v2} {v3};
    \end{tikzpicture}$
    
    As we will calculate later in Example \ref{deletion-contraction}, $X_{\Sigma} = p_{1,0}^3-p_{1,0}p_{1,1}-2p_{1,0}p_{2,0}+2p_{2,1}+p_3-x_0^3$. In the augmented elementary $B$-symmetric basis this is $X_{\Sigma} = (e_1e_2) +(3e_3+q_{1,1}e_1)+(2q_{2,1}+z^3)$ and so $\Sigma$ has 1 acyclic orientation with 2 sinks, $3+1=4$ acyclic orientations with 1 sink and $2+1=3$ acyclic orientations with 0 sinks. Indeed, $\Sigma$ has 8 orientations total.
    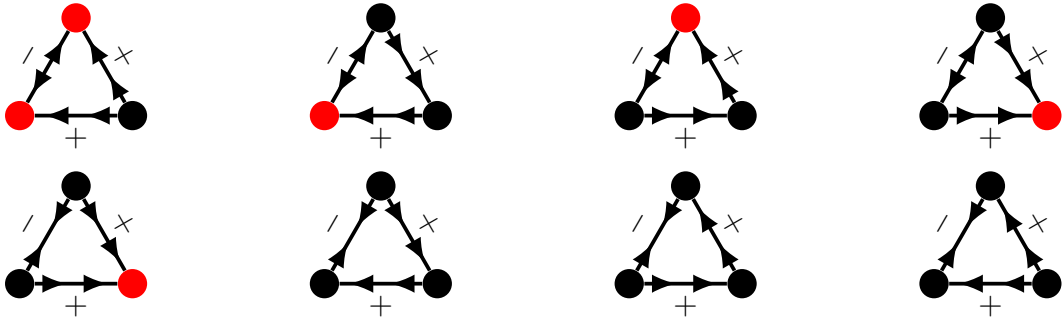
\begin{figure}[h!]\begin{tabular}{c@{\hspace{2cm}}c@{\hspace{2cm}}c@{\hspace{2cm}}c}
    \begin{tikzpicture}
        \node[vertex, red] (v1) at (0,0) {};
        \node[vertex] (v2) at (1.5,0) {};
        \node[vertex, red] (v3) at (0.75, 1.3) {};
        \rtl {v1} {v2} [] [below];
        \ntp {v1} {v3};
        \rtl {v2} {v3};
    \end{tikzpicture} &
    \begin{tikzpicture}
        \node[vertex, red] (v1) at (0,0) {};
        \node[vertex] (v2) at (1.5,0) {};
        \node[vertex] (v3) at (0.75, 1.3) {};
        \rtl {v1} {v2} [] [below];
        \ntp {v1} {v3};
        \ltr {v2} {v3};
    \end{tikzpicture} &
    \begin{tikzpicture}
        \node[vertex] (v1) at (0,0) {};
        \node[vertex] (v2) at (1.5,0) {};
        \node[vertex, red] (v3) at (0.75, 1.3) {};
        \ltr {v1} {v2} [] [below];
        \ntp {v1} {v3};
        \rtl {v2} {v3};
    \end{tikzpicture} &
    \begin{tikzpicture}
        \node[vertex] (v1) at (0,0) {};
        \node[vertex, red] (v2) at (1.5,0) {};
        \node[vertex] (v3) at (0.75, 1.3) {};
        \ltr {v1} {v2} [] [below];
        \ntp {v1} {v3};
        \ltr {v2} {v3};
    \end{tikzpicture} \\
    \begin{tikzpicture}
        \node[vertex] (v1) at (0,0) {};
        \node[vertex, red] (v2) at (1.5,0) {};
        \node[vertex] (v3) at (0.75, 1.3) {};
        \ltr {v1} {v2} [] [below];
        \ptn {v1} {v3};
        \ltr {v2} {v3};
    \end{tikzpicture} &
    \begin{tikzpicture}
        \node[vertex] (v1) at (0,0) {};
        \node[vertex] (v2) at (1.5,0) {};
        \node[vertex] (v3) at (0.75, 1.3) {};
        \rtl {v1} {v2} [] [below];
        \ptn {v1} {v3};
        \ltr {v2} {v3};
    \end{tikzpicture} &
    \begin{tikzpicture}
        \node[vertex] (v1) at (0,0) {};
        \node[vertex] (v2) at (1.5,0) {};
        \node[vertex] (v3) at (0.75, 1.3) {};
        \ltr {v1} {v2} [] [below];
        \ptn {v1} {v3};
        \rtl {v2} {v3};
    \end{tikzpicture} &
    \begin{tikzpicture}
        \node[vertex] (v1) at (0,0) {};
        \node[vertex] (v2) at (1.5,0) {};
        \node[vertex] (v3) at (0.75, 1.3) {};
        \rtl {v1} {v2} [] [below];
        \ptn {v1} {v3};
        \rtl {v2} {v3};
    \end{tikzpicture}
    \end{tabular}
    \caption{The 8 orientations of $\Sigma$ with sinks marked in red.}
    \end{figure}
\end{example}
We may also consider a similar basis to arrive at a similar theorem.
\begin{definition}
    Let $\xi_n = \sum_{a=1}^n\ncr{n}{a}p_{a, 0}$ for $n\geq 1$. 
\end{definition}
Alternatively, we may write this as $p_{n,0} = \sum_{i=1}^n{\ncr{n}{i}}(-1)^{n-i}\xi_i$.

\begin{theorem}\label{main_2}
    If the chromatic $B$-symmetric function of a signed graph $\Sigma$ is written in terms of the basis $\set{\xi_n\st n\geq 1}\cup\set{q_{a,b}\st a,b\geq 1}\cup\{z\}$, then the number of acyclic orientations of $\Sigma$ with $k$ sinks is the sum of the coefficients of terms such that the sum of the indices of each $\xi_n$ factor is equal to $k$.
\end{theorem}
\begin{example}\label{main_thm_2_ex}
    Again, let $\Sigma$ be the signed graph
    \begin{tikzpicture}[baseline = 0.25cm]
        \node[small] (v1) at (0,0) {};
        \node[small] (v2) at (0.75,0) {};
        \node[small] (v3) at (0.375,0.65) {};
        \pstv {v1} {v2} [] [below];
        \ngtv {v1} {v3};
        \pstv {v2} {v3};
    \end{tikzpicture}
    
    We know that $X_{\Sigma} = p_{1,0}^3-p_{1,0}p_{1,1}-2p_{1,0}p_{2,0}+2p_{2,1}+p_3-x_0^3$. In this second basis, this is $X_{\Sigma} = (\xi_1\xi_1\xi_1+\xi_3-2\xi_1\xi_2)+(4\xi_1\xi_2-3\xi_2)+(q_{1,1}\xi_1+3\xi_1)+(2q_{2,1}+z^3)$ and so $\Sigma$ has $1+1-2=0$ acyclic orientations with 3 sinks, $4-3 = 1$ acyclic orientation with 2 sinks, $1+3=4$ acyclic orientations with 1 sink and $2+1=3$ acyclic orientations with 0 sinks.
\end{example}

These theorems specialize to results about unsigned graphs. Specifically, we will show later that Stanley's result about the elementary symmetric basis \cite[Theorem 3.3]{stanley95} follows immediately from Theorem \ref{main_1}. Also, from Theorem \ref{main_2} we obtain the following result 
\begin{definition}
     Let $p_a = \sum_{i\geq 0}x_i^a$ and $\zeta_n = \sum_{a=1}^n{\ncr{n}{a}}p_{a}$, i.e. $p_{a} = \sum_{i=1}^n{\ncr{n}{i}}(-1)^{n-i}\zeta_i$.
\end{definition}
\begin{theorem}
If the chromatic symmetric function of an unsigned graph $G$ is written in terms of the basis $\set{\zeta_n\st n\geq 1}$, then the number of acyclic orientations of $G$ with $k$ sinks is the sum of the coefficients of terms such that the sum of the indices of each $\zeta_n$ factor is equal to $k$.
\end{theorem}

\begin{example}\label{main_thm_3_ex}
    Let $G$ be the unsigned graph $\begin{tikzpicture}[baseline = 0.125cm]
        \node[small] (v1) at (0,0) {};
        \node[small] (v2) at (0.75,0) {};
        \node[small] (v3) at (0.375,0.65) {};
        \edge {v1} {v2} [] [below];
        \edge {v1} {v3};
        \edge {v2} {v3};
    \end{tikzpicture}$
    
    Then 
$X_G = p_1^3-3p_1p_2+2p_3$. This can also be written as 
\begin{align*} 
X_G =& \zeta_1^3-3\zeta_1(\zeta_2-2\zeta_1)+2(\zeta_3-3\zeta_2+3\zeta_1)\\
& = (\zeta_1\zeta_1\zeta_1-3\zeta_1\zeta_2+2\zeta_3)+(6\zeta_1\zeta_1-6\zeta_2)+(6\zeta_1)
\end{align*}
 So $G$ has $1-3+2=0$ acyclic orientations with 3 sinks, $6-6=0$ acyclic orientations with 2 sinks and 6 acyclic orientations with 1 sink. 
\end{example}

Lastly, Theorem \ref{main_2} also gives an equivalent form of Zaslavsky's result \cite [Corollary 4.1]{zaslavsky}. The sum of the absolute value of the coefficients of $X_\Sigma$ written in the $p$-basis is equal to the total number of acyclic orientations of $\Sigma$.

\begin{example}\label{main_thm_4_ex}
    Let $\Sigma$ be the signed graph$\begin{tikzpicture}[baseline = 0.25cm]
        \node[small] (v1) at (0,0) {};
        \node[small] (v2) at (0.75,0) {};
        \node[small] (v3) at (0.375,0.65) {};
        \pstv {v1} {v2} [] [below];
        \ngtv {v1} {v3};
        \pstv {v2} {v3};
    \end{tikzpicture}$
    
    Then 
$$
        X_{\Sigma} = p_{1,0}^3-p_{1,0}p_{1,1}-2p_{1,0}p_{2,0}+2p_{2,1}+p_{3,0}-x_0^3
$$
    So $\Sigma$ has $|1|+|-1|+|-2|+|2|+|1|+|-1| = 8$ acyclic orientations.
\end{example}
We can also see that this holds for unsigned graphs as well

\begin{example}\label{main_thm_5_ex}
    Let $G$ be the unsigned graph $\begin{tikzpicture}[baseline = 0.25cm]
        \node[small] (v1) at (0,0) {};
        \node[small] (v2) at (0.75,0) {};
        \node[small] (v3) at (0.375,0.65) {};
        \edge {v1} {v2} [] [below];
        \edge {v1} {v3};
        \edge {v2} {v3};
    \end{tikzpicture}$
    
    Then 
\[
    X_G = p_1^3-3p_1p_2+2p_3.
\]
 So, $G$ has $|1|+|-3|+|2|= 6$ acyclic orientations. 
\end{example}

%% file: Sections/deletion_contraction.tex
\section{Weighted Deletion-Contraction}
Now we will take a moment to consider how we might be able to calculate $X_{\Sigma}$. 

The classical chromatic polynomial satisfies a relation known as ``deletion-contraction'' \cite{stanley73}, which turns out to be immensely useful for both computation and theoretical considerations.
Stanley's chromatic symmetric function satisfies a similar relation, but this time involving weighted graphs. This work appears in \cite{Chmutov_paper, Noble} and was recently rediscovered in \cite{Crew}.
Due to unpublished work of Enouen, Fawcett, Raghavan, and Shah \cite{Chmutov_slides} we have a generalization of this for the chromatic $B$-symmetric function, namely a ``doubly weighted'' deletion-contraction rule for signed graphs.

For a signed graph $\Sigma$, a \emph{double weight function}, $w = (w_+,w_-)$, will mean a pair of functions $w_+,w_-: V(\Sigma) \to \N$.
For the doubly weighted signed graph $(\Sigma, w)$, we define
\[
    X_{(\Sigma, w)} = \sum_{\kappa \in \mathcal{P}(G)}x_{\kappa(v_1)}^{w_+(v_1)} x_{-\kappa(v_1)}^{w_-(v_1)}\dots x_{\kappa(v_n)}^{w_+(v_n)} x_{-\kappa(v_n)}^{w_-(v_n)}.
\]
We observe that this extends the chromatic $B$-symmetric function defined above if we treat an unweighted signed graph $\Sigma$ as having a double weight function satisfying $w(v) = (1,0)$ for all $v\in V(\Sigma)$.

To define the weighted deletion-contraction rule, we must have a notion of both deletion and contraction on signed graphs.
The latter is easy: given a (possibly doubly weighted) signed graph $\Sigma$ and an edge $e$ in $\Sigma$, we write $\Sigma\setminus e$ for the graph obtained by deleting the edge $e$.
If $e \epts uv$ is a positive edge, the \emph{contraction} of $\Sigma$ along $e$ is the graph $\Sigma/e$ whose vertex set is $V(\Sigma)$ modulo the relation $u\sim v$, and whose edges are obtained from $E(\Sigma)\setminus\{e\}$ by replacing all endpoints $u$ and $v$ by the equivalence class $\{u,v\}$.
If $\Sigma$ has a double weight function $w$, then we consider the induced double weight $\widetilde{w}$ on $G/e$, which only differs from $w$ by $\widetilde{w}(\{u,v\}) =  w(u) + w(v)$ (where the addition of pairs is component-wise).

% [examples of deletion and contraction]

% Now that we have defined deletion and contraction, the deletion-contraction rule is the following.

\begin{theorem}\label{thm:dc}
    Let $(\Sigma,w)$ be a doubly weighted signed graph, and suppose $e_0\in E(G)$ is a positive edge, i.e., $\sgn(e_0) = +$.
    Then $X_{(\Sigma,w)} = X_{(\Sigma\setminus e_0,w)} - X_{(\Sigma/e_0,\widetilde{w})}$.
\end{theorem}

\begin{proof}
We write the chromatic $B$-symmetric function as a sum over all (not necessarily proper) colorings $\kappa : V(\Sigma) \to \Z$ as follows:
\[
    X_{(\Sigma,w)} = \sum_{\kappa: V(\Sigma)\to\Z}\left(\prod_{v\in V(\Sigma)} x_{\kappa(v)}^{w_+(v)} x_{-\kappa(v)}^{w_-(v)} \prod_{\substack{e\in E(\Sigma) \\ e\epts uv}} \left(1 - \delta_{\kappa(u)}^{\sgn(e)\kappa(v)}\right)\right)
\]
where $\delta_a^b$ is the Kronecker delta function, $\delta_{a}^b = 1$ if $a=b$ and $\delta_{a}^b=0$ otherwise.
Write $e_0 \epts u_0 v_0$ and expand across the factor $1-\delta_{\kappa(u_0)}^{\kappa(v_0)}$ to get
\begin{align*}
    X_{(\Sigma,w)} &= \sum_{\kappa: V(\Sigma)\to\Z}\left(
    \prod_{v\in V(\Sigma)}x_{\kappa(v)}^{w_+(v)} x_{-\kappa(v)}^{w_-(v)}\prod_{\substack{e\in E(\Sigma)\setminus\{e_0\} \\ e\epts uv}}\left(
    1 - \delta_{\kappa(u)}^{\sgn(e)\kappa(v)}\right)\right) \\
    &\hspace{4em}-\sum_{\substack{\kappa: \Sigma\to\Z \\ \kappa(u) = \kappa(v)}}\left(
    \prod_{v\in V(\Sigma)}
    x_{\kappa(v)}^{w_+(v)} x_{-\kappa(v)}^{w_-(v)}
    \prod_{\substack{e\in E(\Sigma)\setminus\{e_0\} \\ e\epts uv}}
    \left(1 - \delta_{\kappa(u)}^{\sgn(e)\kappa(v)}\right)\right).
\end{align*}
But this is just the deletion-contraction rule we wanted to prove.
\end{proof}

This theorem initially appears to be lacking in usefulness since it can only deal with positive edges.
Instead of expanding the rule to consider negative edges, which would yield too complicated a result, we work around the problem by introducing a method of turning negative edges into positive edges under which the chromatic $B$-symmetric function is invariant.
Given a doubly weighted signed graph $(\Sigma,w)$ and a vertex $v\in V(\Sigma)$, the graph obtained from \emph{switching} at $v$ is the graph $(\Sigma^v,w^v)$, where $w^v$ differs only from $w$ in that if $w(v) = (a,b)$ then $w^v(v) = (b,a)$, and $\Sigma^v$ is the signed graph $\Sigma$ except all positive non-loop edges connected to $v$ are now negative and all negative non-loop edges connected to $v$ are now positive.

\begin{figure}[h!t]
\includegraphics{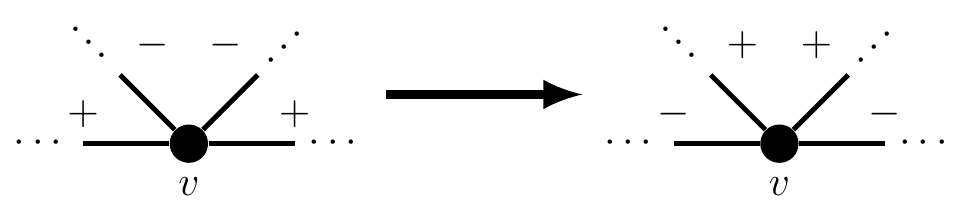}
\caption{Switching a vertex $v$.}
% \begin{tikzpicture}
%     \node[vertex, label=below:$v$] (v1) at (0,0) {};
%     \node (dots1) at (-1.5,0) {$\cdots$};
%     \node[rotate=-45] (dots2) at (-1,1) {$\cdots$};
%     \node[rotate=45] (dots3) at (1,1) {$\cdots$};
%     \node (dots4) at (1.5,0) {$\cdots$};
            
%     \edge {v1} {dots1} node [at end, above] {$+$};
%     \edge {v1} {dots2} node [at end, above right] {$-$};
%     \edge {v1} {dots3} node [at end, above left] {$-$};
%     \edge {v1} {dots4} node [at end, above] {$+$};
            
%     \draw[line width = 1mm, -latex] (2,0.5) -- (4,0.5);
    
%     \node[vertex, label=below:$v$] (v1) at (6,0) {};
%     \node (dots5) at (4.5,0) {$\cdots$};
%     \node[rotate=-45] (dots6) at (5,1) {$\cdots$};
%     \node[rotate=45] (dots7) at (7,1) {$\cdots$};
%     \node (dots8) at (7.5,0) {$\cdots$};
            
%     \edge {v1} {dots5} node [at end, above] {$-$};
%     \edge {v1} {dots6} node [at end, above right] {$+$};
%     \edge {v1} {dots7} node [at end, above left] {$+$};
%     \edge {v1} {dots8} node [at end, above] {$-$};
% \end{tikzpicture}
\end{figure}
We now justify our claim above that the chromatic $B$-symmetric function is invariant under switching.

\begin{lemma}\label{lem:switching}
    Let $(\Sigma,w)$ be a doubly weighted signed graph, and let $v_0 \in V(\Sigma)$.
    Write $(\Sigma^{v_0},w^{v_0})$ for the graph obtained by switching at $(\Sigma,w)$ at $v_0$.
    Then $X_{(\Sigma,w)} = X_{(\Sigma^{v_0},w^{v_0})}$.
\end{lemma}

\begin{proof}
    We write
    \begin{align*}
        X_{(\Sigma^{v_0},w^{v_0})} = \sum_{\kappa: \Sigma \to \Z} & \left(\prod_{v\in V(\Sigma)\setminus\{v_0\}}x_{\kappa(v)}^{w_+(v)}x_{-\kappa(v)}^{w_-(v)}\right)
        \left(x_{-\kappa(v_0)}^{w_+(v_0)} x_{\kappa(v_0)}^{w_-(v_0)}\right) \\
        &\cdot\prod_{\substack{e \in E(G) \\ e\epts uv\\ v_0 \notin \{u,v\}}}\left(1-\delta_{\kappa(u)}^{\sgn(e)\kappa(v)}\right)\prod_{\substack{e \in E(G)\\ e \epts uv_0}}\left(1-\delta_{\kappa(u)}^{-\sgn(e)\kappa(v_0)}\right).
    \end{align*}
    If we reindex the sum, by replacing each $\kappa$ with the function $\kappa^{v_0}: \Sigma \to\Z$ differing from $\kappa$ only by $\kappa^{v_0}(v_0) = -\kappa(v_0)$, we see that the whole expression is just $X_G$.
\end{proof}

Now we can perform deletion-contraction on a graph, switching negative edges to positive when necessary, until we end up with disjoint vertices, some of which may have a negative loop. The chromatic $B$-symmetric function of disjoint vertices is the product of the chromatic $B$-symmetric function of each of the individual vertices and if $(\Sigma,w)$ is a doubly weighted signed graph which consists of a single vertex of weight $(a,b)$, then $X_{(\Sigma,w)} = \sum_{i\in \Z}x_{i}^ax_{-i}^b=p_{a,b}$. If this vertex has a negative loop then we have 
\[
    X_{(\Sigma,w)} = \sum_{i\in \Z\setminus \{0\}}x_{i}^ax_{-i}^b = \left(\sum_{i\in \Z}x_{i}^ax_{-i}^b\right)-x_0^{a+b} = p_{a,b}-x_0^{a+b}.
\]
% We will use the notation $p_{a,b} \coloneqq \sum_{i\in \Z}x_{i}^ax_{-i}^b$ and we will refer to the set $\set{p_{a,b}\st a\geq b\geq 0}\cup\set{x_0}$ as the \emph{power basis}.
It follows that $X_{(\Sigma,w)}$ and hence $X_{\Sigma}$ can always be written in terms of the set $\set{p_{a,b}\st a\geq b\geq 0}\cup\set{x_0}$, which forms the \emph{power basis}.

% We've seen that the $B$-symmetric chromatic function of a single vertex with no loops and weight $(a,b)$ has the form $p_{a,b}=\sum_{i\in\Z}x_i^ax_{-1}^b$, and the $B$-symmetric chromatic function of a single vertex with a negative loop and weight $(a,b)$ has the form $p_{a,b}-x_0^{a+b}=\sum_{i\in\Z}x_i^ax_{-i}^b-x_0^{a+b}$. (A single vertex with positive loop is assigned the zero function.) Note that $p_{a,b}=p_{b,a}$, so without loss of generality we may consider $a\geq b$. These $p$ functions are enough to form a basis of the space of $B$-symmetric functions: we call $\set{p_{a,b}\st a\geq b\geq0}\cup\set{x_0}$ the \emph{power basis}.

\begin{example}\label{deletion-contraction}
    Take $\Sigma$ to be the triangle with one negative edge and for convenience we will write just a doubly weighted signed graph in place of writing the graph and its double weight function in the subscript of $X$.

    Then we have
    \begin{align*}
        \begin{tikzpicture}[baseline = 0.35cm]
            \node[small, label=below left:{$(1,0)$}] (v1) at (0,0) {};
            \node[small, label=below right:{$(1,0)$}] (v2) at (1,0) {};
            \node[small, label=above:{$(1,0)$}] (v3) at (0.5,0.866) {};
            \pstv {v1} {v2} [] [below];
            \pstv {v1} {v3};
            \ngtv {v2} {v3};
        \end{tikzpicture}\hspace{-14pt}& \\=&
        \begin{tikzpicture}[baseline = 0.35cm]
            \node[small, label=below:{$(1,0)$}] (v1) at (0,0) {};
            \node[small, label=below:{$(1,0)$}] (v2) at (1,0) {};
            \node[small, label=above:{$(1,0)$}] (v3) at (0.5,0.866) {};
            \pstv {v1} {v3};
            \ngtv {v2} {v3};
        \end{tikzpicture}
        -
        \begin{tikzpicture}[baseline = 0.35cm]
            \node[small, label=below:{$(2,0)$}] (v1) at (0.5,0) {};
            \node[small, label=above:{$(1,0)$}] (v3) at (0.5,0.866) {};
            \draw[bend left = 60, line width = 0.5mm] (v1) to (v3) node[midway, label = above: $+$] {};
            \draw[bend right = 60, line width = 0.5mm] (v1) to (v3);
            \node[rotate = 90] (sign) at (1,0.45) {$-$};
        \end{tikzpicture}
        \\=
        &\begin{tikzpicture}[baseline = 0.35cm]
            \node[small, label=below:{$(1,0)$}] (v1) at (0,0) {};
            \node[small, label=below:{$(1,0)$}] (v2) at (1,0) {};
            \node[small, label=above:{$(1,0)$}] (v3) at (0.5,0.866) {};
            \ngtv {v2} {v3};
        \end{tikzpicture}
        -
        \begin{tikzpicture}[baseline = 0.35cm]
            \node[small, label=below:{$(1,0)$}] (v2) at (1,0) {};
            \node[small, label=above:{$(2,0)$}] (v3) at (0.5,0.866) {};
            \ngtv {v2} {v3};
        \end{tikzpicture}
        -
        \begin{tikzpicture}[baseline = 0.35cm]
            \node[small, label=below:{$(2,0)$}] (v1) at (0.5,0) {};
            \node[small, label=above:{$(1,0)$}] (v3) at (0.5,0.866) {};
            \ngtv {v1} {v3} [] [sloped];
            \node[rotate = 90] (sign) at (1,0.45) {$-$};
        \end{tikzpicture}
        +
        \begin{tikzpicture}[baseline = 0cm]
            \node[small, label=below:{$(3,0)$}] (v1) at (0,0) {};
            \draw[line width=0.5mm, loop above, out = 60, in =120, looseness = 15] (v1) to (v1);
            \node (sign) at (0,0.75) {$-$};
        \end{tikzpicture}
        \\
        = &\begin{tikzpicture}[baseline = 0.35cm]
            \node[small, label=below:{$(1,0)$}] (v1) at (0,0) {};
            \node[small, label=below:{$(1,0)$}] (v2) at (1,0) {};
            \node[small, label=above:{$(1,0)$}] (v3) at (0.5,0.866) {};
        \end{tikzpicture}
        - 
        \begin{tikzpicture}[baseline = 0.35cm]
            \node[small, label=below:{$(1,0)$}] (v1) at (0,0) {};
            \node[small, label=above:{$(1,1)$}] (v3) at (0.5,0.866) {};
        \end{tikzpicture}
        -
        \begin{tikzpicture}[baseline = 0.35cm]
            \node[small, label=below:{$(1,0)$}] (v2) at (1,0) {};
            \node[small, label=above:{$(2,0)$}] (v3) at (0.5,0.866) {};
        \end{tikzpicture}
        +
        \begin{tikzpicture}[baseline = 0cm]
            \node[small, label=below:{$(2,1)$}] (v2) at (1,0) {};
        \end{tikzpicture}
        -
        \begin{tikzpicture}[baseline = 0.35cm]
            \node[small, label=below:{$(2,0)$}] (v1) at (0.5,0) {};
            \node[small, label=above:{$(1,0)$}] (v3) at (0.5,0.866) {};
        \end{tikzpicture}
        +
        \begin{tikzpicture}[baseline = 0cm]
            \node[small, label=below:{$(2,1)$}] (v1) at (0.5,0) {};
        \end{tikzpicture}
        +
        \begin{tikzpicture}[baseline = 0cm]
            \node[small, label=below:{$(3,0)$}] (v1) at (0,0) {};
            \draw[line width=0.5mm, loop above, out = 60, in =120, looseness = 15] (v1) to (v1);
            \node (sign) at (0,0.75) {$-$};
        \end{tikzpicture}
        \end{align*}
        
        And so 
    \begin{align*}
        X_{\Sigma} 
        &=  p_{1,0}^3 - p_{1,0}p_{1,1} - p_{1,0}p_{2,0}+p_{2,1}-p_{1,0}p_{2,0}+p_{2,1}+p_{3,0}-x_0^3\\
        &= p_{1,0}^3 -p_{1,0}p_{1,1}- 2p_{1,0}p_{2,0}+2p_{2,1}+p_{3,0}-x_0^3
    \end{align*}
\end{example}

%% file: Sections/covering.tex
\subsection{The Covering Graph}
In working with signed graphs, it is exceptionally helpful to consider a construction due to Zaslavsky called the \emph{covering graph} (called ``\emph{signed} covering graph" in \cite[Theorem 3.1]{zaslavsky}).
%\footnote{Zaslavsky introduced the object as the \emph{signed} covering graph, but we feel the term is confusing since the covering graph is unsigned. As such, we drop the word ``signed''.}.
Given a signed graph $\Sigma$, the covering graph of $\Sigma$, denoted $\overline{\Sigma}$, is an unsigned graph such that for every vertex $v$ in $\Sigma$, there are vertices $+v$ and $-v$ in $\overline{\Sigma}$. Also, for every edge $e \epts uv$ in $\Sigma$, there are two edges, one between $+v$ and $\sgn(e)u$ and another one between $-v$ and $-\sgn(e)u$ in $\overline{\Sigma}$.

Additionally, given an orientation on $\Sigma$, we may induce an orientation on $\overline{\Sigma}$ such that for every arrow incident with a vertex $v$ in $\Sigma$, the corresponding arrow is incident with $+v$ in $\overline{\Sigma}$ and the reverse arrow is incident with $-v$. 
\begin{figure}[h!t]
    \includegraphics{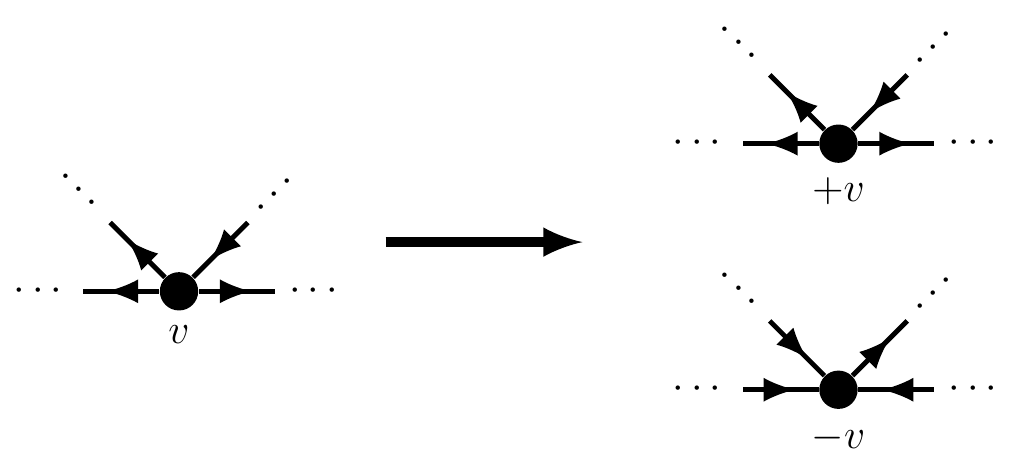}
    \caption{Creating an oriented covering graph.}
    % \begin{tikzpicture}
    %     \node[vertex, label=below:$v$] (v1) at (-0.1,0) {};
    %     \node (dots1) at (-1.5,0) {$\cdots$};
    %     \node[rotate=-45] (dots2) at (-1.1,1) {$\cdots$};
    %     \node[rotate=45] (dots3) at (0.9,1) {$\cdots$};
    %     \node (dots4) at (1.3,0) {$\cdots$};
        
    %     \drct {v1} {dots1} {latex-} node [at end, above] {};
    %     \drct {v1} {dots2} {latex-} node [at end, above right] {};
    %     \drct {v1} {dots3} {latex-} node [at end, above left] {};
    %     \drct {v1} {dots4} {-latex} node [at end, above] {};
        
    %     \draw[line width = 1mm, -latex] (2,0.5) -- (4,0.5);
        
    %     \node[vertex, label=below:$+v$] (v3) at (6.6,1.5) {};
    %     \node (dots1) at (5.2,1.5) {$\cdots$};
    %     \node[rotate=-45] (dots2) at (5.6,2.5) {$\cdots$};
    %     \node[rotate=45] (dots3) at (7.6,2.5) {$\cdots$};
    %     \node (dots4) at (8,1.5) {$\cdots$};
        
    %     \drct {v3} {dots1} {latex-} node [at end, above] {};
    %     \drct {v3} {dots2} {latex-} node [at end, above right] {};
    %     \drct {v3} {dots3} {latex-} node [at end, above left] {};
    %     \drct {v3} {dots4} {-latex} node [at end, above] {};
        
    %     \node[vertex, label=below:$-v$] (v2) at (6.6,-1) {};
    %     \node (dots5) at (5.2,-1) {$\cdots$};
    %     \node[rotate=-45] (dots6) at (5.6,0) {$\cdots$};
    %     \node[rotate=45] (dots7) at (7.6,0) {$\cdots$};
    %     \node (dots8) at (8,-1) {$\cdots$};
    
    %     \drct {v2} {dots5} {-latex} node [at end, above] {};
    %     \drct {v2} {dots6} {-latex} node [at end, above right] {};
    %     \drct {v2} {dots7} {-latex} node [at end, above left] {};
    %     \drct {v2} {dots8} {latex-} node [at end, above] {};
    % \end{tikzpicture}
\end{figure}

As it turns out, on any given edge in $\overline{\Sigma}$ and for any orientation induced this way, both arrows on the edge will be pointing in the same direction and so we will replace the two arrows on an edge with a single arrow pointing in their shared direction.

\begin{figure}[!ht]
    \includegraphics{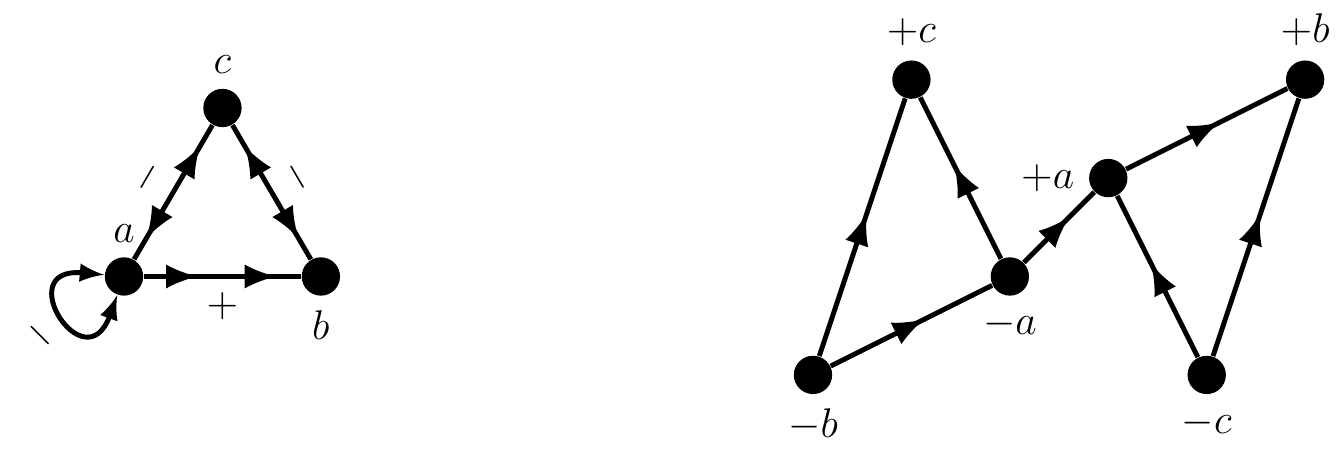}
    \caption{An oriented signed graph $\Sigma$ on the left and its covering graph $\overline{\Sigma}$ with the induced orientation on the right.}
    % \begin{tikzpicture}
    %     \node[vertex, label = above:$a$] (v1) at (0,0) {};
    %     \node[vertex, label = below:$b$] (v2) at (2,0) {};
    %     \node[vertex, label = above:$c$] (v3) at (1,1.712) {};
    
    %     \ltr {v1} {v2} [] [below];
    %     \ntp {v1} {v3} [] [above] [sloped];
    %     \ptn {v2} {v3} [] [above] [sloped];
    %     \draw[latex-latex, line width=0.5mm, loop below, out = 175, in = 250, looseness = 10] (v1) to (v1);
    %     \node[rotate = -45] (sign) at (-0.85,-0.6) {$-$};
    
    %     % \draw[line width = 1mm, -latex] (3.5,1) -- (5.5,1);
        
    %     \node[vertex, label = above:$+c$] (+c) at (8,2) {};
    %     \node[vertex, label = above:$+b$] (+b) at (12,2) {};
    %     \node[vertex, label = left:$+a$] (+a) at (10,1) {};
    %     \node[vertex, label = below:$-c$] (-c) at (11,-1) {};
    %     \node[vertex, label = below:$-a$] (-a) at (9,0) {};
    %     \node[vertex, label = below:$-b$] (-b) at (7,-1) {};
    %     \drct {-b} {+c} {-latex};
    %     \drct {-a} {+c} {latex-};
    %     \drct {-a} {+a} {-latex};
    %     \drct {-c} {+a} {latex-};
    %     \drct {+a} {+b} {-latex};
    %     \drct {-b} {-a} {-latex};
    %     \drct {-c} {+b} {-latex};
    % \end{tikzpicture}
\end{figure}

% Next we will consider the B-Symmetric Chromatic Function, $X_{\Sigma}$ of a signed graph $\Sigma$, as introduced by \cite{Egge}. In analogy with the unsigned case, we put
% \[
%     X_\Sigma(\dots,x_{-1},x_0,x_1,\dots) = \sum_{\kappa\in\mathcal P(\Sigma)}x_{\kappa(v_1)}x_{\kappa(v_2)}\cdots x_{\kappa(v_n)}
% \]
% where $v_1,v_2,\dots, v_n$ are the vertices of $\Sigma$ and $\mathcal{P}(\Sigma)$ is the set of proper colorings of $\Sigma$. For notational convenience, we will put $x^{\kappa} = x_{\kappa(v_1)}x_{\kappa(v_2)}\cdots x_{\kappa(v_n)}$ whenever applicable. 

%% file: Sections/linear_extensions.tex
\subsection{Orientation preserving colorings}
%Now that we have a method for calculating $B$-symmetric chromatic functions, we move on to more theoretical considerations.
Following Zaslavsky's lead \cite{zaslavsky}, we will consider the connection between colorings and orientations. Each proper coloring of a signed graph $\Sigma$ induces (or \emph{preserves}) a unique acyclic orientation and moreover, each acyclic orientation can be induced from some proper coloring. Given a proper coloring of $\Sigma$ we can construct an acyclic orientation of $\Sigma$ by first constructing a coloring on the covering graph $\overline{\Sigma}$ such that for every vertex $v$ in $\Sigma$, the vertex $+v$ in $\overline{\Sigma}$ has the same color as $v$ and the vertex $-v$ in $\overline{\Sigma}$ has the negative of the color of $v$. With this new coloring on $\overline{\Sigma}$, we create an orientation on $\overline{\Sigma}$ such that lower colors always point to higher colors and we use this to give us an orientation on $\Sigma$. Proof that an orientation constructed this way is acyclic is given in \cite{zaslavsky}.

\begin{figure}[h!]
    \noindent
    \includegraphics{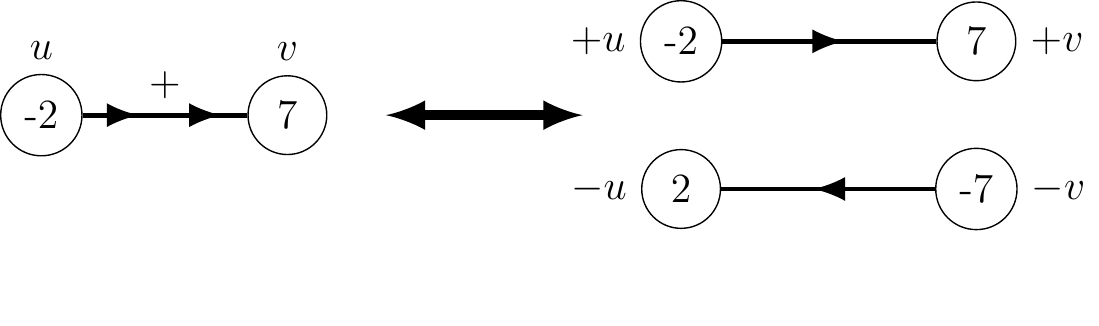}
    % \begin{tikzpicture}
    %     \node[colored, label=above:$u$] (u) at (0,0) {-2};
    %     \node[colored, label=above:$v$] (v) at (2.5,0) {$7$};
    %     \draw[line width = 1mm, latex-latex] (3.5,0) -- (5.5,0);
    %     \node[colored, label=left:$+u$] (+u) at (6.5,0.75) {-2};
    %     \node[colored, label=right:$+v$] (+v) at (9.5,0.75) {7};
    %     \node[colored, label=left:$-u$] (-u) at (6.5,-0.75) {2};
    %     \node[colored, label=right:$-v$] (-v) at (9.5,-0.75)
    %     {-7};
        
    %     \ltr {u} {v} [] [above];
    %     \drct {+u} {+v} {-latex};
    %     \drct {-u} {-v} {latex-};
        
    %     \node (spurious) at (0,-2) {};
    %     %\draw[dashed, line width = 0.5mm] (-0.5,-1.4) -- (10,-1.4);
    % \end{tikzpicture}
    \noindent
    \includegraphics{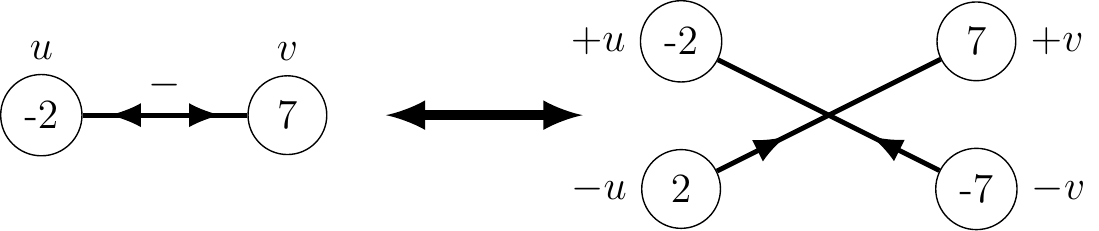}
    \caption{Two examples of colorings preserving orientations
    % Vertex $u$ is colored $-2$ and vertex $v$ is colored $7$. If $u$ and $v$ are connected by a positive edge, then that edge must point from $u$ to $v$ so that on the covering graph, $+u$ points to $+v$ and $-v$ points to $-u$, since $+u$ has a lower color than $+v$ and $-v$ has a lower color than $-u$. Similarly, if $u$ and $v$ are connected by a negative edge, then the arrows must point outward so that on the covering graph, $-u$ points to $+v$ and $-v$ points to $+u$, since $-u$ has a lower color than $+v$ and $-v$ has a lower color than $+u$.}
    }
    % \begin{tikzpicture}
    %     \node[colored, label=above:$u$] (u) at (0,0) {-2};
    %     \node[colored, label=above:$v$] (v) at (2.5,0) {$7$};
        
    %     \draw[line width = 1mm, latex-latex] (3.5,0) -- (5.5,0);
        
    %     \node[colored, label=left:$+u$] (+u) at (6.5,0.75) {-2};
    %     \node[colored, label=right:$+v$] (+v) at (9.5,0.75) {7};
    %     \node[colored, label=left:$-u$] (-u) at (6.5,-0.75) {2};
    %     \node[colored, label=right:$-v$] (-v) at (9.5,-0.75) {-7};
        
    %     \ntp {u} {v} [] [above];
    %     \drct {+u} {-v} {latex-} [] [near end];
    %     \drct {-u} {+v} {-latex} [] [near start];
    % \end{tikzpicture}
\end{figure}

For an orientation $P$ of a signed graph $\Sigma$, we will we use $\mathcal{C}(P)$ to denote the set of all proper colorings of $\Sigma$ which preserve $P$. Additionally, we will use $\mathcal{A}(\Sigma)$ to denote the set of acyclic orientations of $\Sigma$.

Now it is easy to see that 
$$
\mathcal{P}(\Sigma) = \bigsqcup_{P\in \mathcal{A}(\Sigma)}\mathcal{C}(P)
$$
In other words, we can partition the set of all proper colorings of $\Sigma$ into sets corresponding to the acyclic orientations of $\Sigma$.

If, for an acyclic orientation $P$ of $\Sigma$, we put $Y_{P} = \sum_{\kappa\in \mathcal{C}(P)}x^{\kappa}$, then this gives us that 
$$
X_{\Sigma} = \sum_{\kappa\in \mathcal{P}(\Sigma)}x^{\kappa} = \sum_{P\in \mathcal{A}(\Sigma)}\left(\sum_{\kappa \in\mathcal{C}(P)}x^{\kappa}\right) = \sum_{P\in\mathcal{A}(\Sigma)}Y_P
$$
\subsection{Linear extensions of signed posets}
It will be very useful to partition the proper colorings of $\Sigma$ even further. In Stanley's treatment of the unsigned case \cite{stanley95} he does this by treating acyclic oriented graphs as posets of their vertices and considering \emph{linear extensions} (total orderings which contain all relations in the poset) of these posets. In analogy with this, we will consider a linear extension of an orientation of a signed graph to be a linear extension (in the unsigned sense) of its oriented covering graph such that if a vertex $+v$ is $k$th from the top of the total order, then $-v$ is $k$th from the bottom of the total order. 
\begin{figure}[b]
%    \includegraphics{tikz/Linear_Extension.pdf}
%     \caption{A signed graph, its corresponding covering graph, and one possible linear extension of the covering graph.}
% \end{figure}

 \centering
    \label{fig:linear_extension}

\begin{tikzpicture}
    \node[vertex, label = above:$a$] (v1) at (0,0) {};
    \node[vertex, label = below:$b$] (v2) at (2,0) {};
    \node[vertex, label = above:$c$] (v3) at (1,1.712) {};

    \ltr {v1} {v2} [] [below];
    \ntp {v1} {v3} [] [above] [sloped];
    \ptn {v2} {v3} [] [above] [sloped];
    \draw[latex-latex, line width=0.5mm, loop below, out = 175, in = 250, looseness = 10] (v1) to (v1);
    \node[rotate = -45] (sign) at (-0.85,-0.6) {$-$};

    \draw[line width = 1mm, -latex] (3,1) -- (5,1) node[midway, label = above: Covering Graph] {};
    
    \node[vertex, label = above:$+c$] (+c) at (6.2,2) {};
        \node[vertex, label = above:$+b$] (+b) at (10.2,2) {};
        \node[vertex, label = left:$+a$] (+a) at (8.2,1) {};
        \node[vertex, label = below:$-c$] (-c) at (9.2,-1) {};
        \node[vertex, label = below:$-a$] (-a) at (7.2,0) {};
        \node[vertex, label = below:$-b$] (-b) at (5.2,-1) {};
        \drct {-b} {+c} {-latex};
        \drct {-a} {+c} {latex-};
        \drct {-a} {+a} {-latex};
        \drct {-c} {+a} {latex-};
        \drct {+a} {+b} {-latex};
        \drct {-b} {-a} {-latex};
        \drct {-c} {+b} {-latex};
    
    \draw[line width = 1mm, -latex] (11,1) -- (13,1) node[midway, label = above: Linear Extension] {};
    
    \node[vertex, label = right:$-c$] (+a1) at (14,1.5) {};
        \node[vertex, label = right:$+a$] (+b1) at (14,2.5) {};
        \node[vertex, label = right:$+b$] (+c1) at (14,3.5) {};
        \node[vertex, label = right:$+c$] (-a1) at (14,0.5) {};
        \node[vertex, label = right:$-a$] (-b1) at (14,-0.5) {};
        \node[vertex, label = right:$-b$] (-c1) at (14,-1.5) {};
        \drct {+a1} {+b1} {-latex};
        \drct {+b1} {+c1} {-latex};
        \drct {-a1} {+a1} {-latex};
        \drct {-b1} {-a1} {-latex};
        \drct {-c1} {-b1} {-latex};
\end{tikzpicture}
\caption{A signed graph, its corresponding covering graph, and one possible linear extension of the covering graph.}
\end{figure}
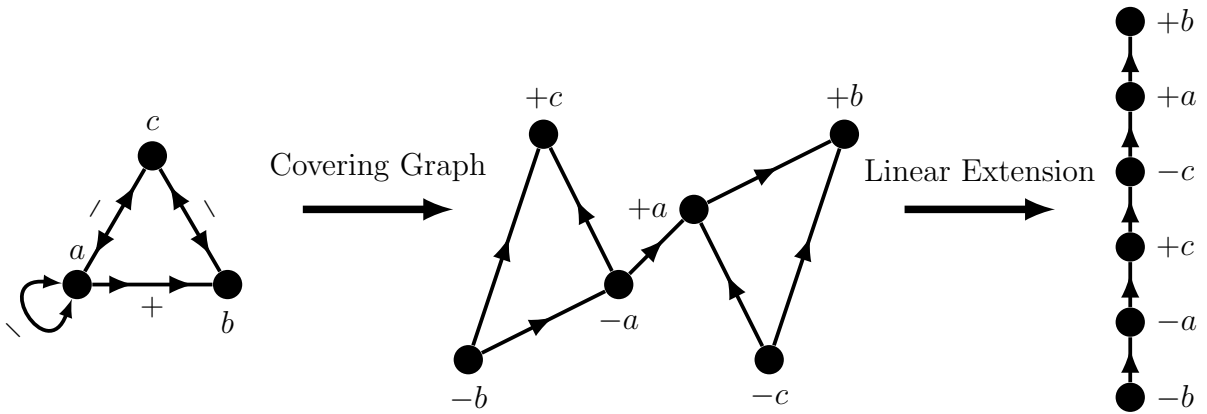
If $P$ is an acyclic orientation of some signed graph, then we will call $P$ a $\emph{signed poset}$. For a signed poset $P$, we will use the notation $\epsilon_1 v <_P \epsilon_2 u$ to mean that $\epsilon_1 v$ points to $\epsilon_2 u$ in the oriented covering graph, where $\epsilon_1, \epsilon_2\in \{+,-\}$. We will retain the convention $\epsilon_1, \epsilon_2\in \{+,-\}$ for the rest of the paper.

We will find it very convenient to treat linear extensions as functions. For a signed poset $P$ on a graph $\Sigma$, which has $d$ vertices, we can define a linear extension of $P$ to be a function $\alpha$ from the vertex set of  $\overline{\Sigma}$ to the set $\{-d,-d+1,\dots, -1, 1, \dots, d-1,d\}$ such that $\epsilon_1 u >_P \epsilon_2 v \implies \alpha(\epsilon_1 u) > \alpha( \epsilon_2 v)$. Additionally we require that $\alpha(\epsilon v) = \epsilon\alpha(+v)$ for all $v$. We will also put $\sgn_\alpha(v)=\frac{|\alpha(+v)|}{\alpha(+v)}$.

\begin{definition}
    For a signed poset $P$, let $\mathcal{L}(P)$ denote the set of all linear extensions of $P$.
\end{definition}
\begin{definition}
     For a proper coloring $\kappa$ and a vertex $v$, let $\sgn_{\kappa}(v) = \frac{|\kappa(v)|}{\kappa(v)}$ when $\kappa(v)\neq 0$.
\end{definition}
\begin{definition}
    Given a signed poset $P$ and $\alpha,\omega\in \mathcal{L}(P)$, define $\mathcal K(\alpha,\omega)$ to be the set of all colorings $\kappa$ such that for all vertices $u,v$:
    \begin{enumerate}[label={(4.6.\arabic*)}]
            \item If $|\alpha(\epsilon_1u)| < |\alpha(\epsilon_2v)|$, then $|\kappa(u)| \leq |\kappa(v)|$.
            \item $\sgn_\alpha(v) = \sgn_\kappa(v)$ or $\kappa(v)=0$.
            \item If both $\alpha(\epsilon_1 u) < \alpha(\epsilon_2 v)$ and $\omega(\epsilon_1 u) < \omega(\epsilon_2 v)$ are true, then $\epsilon_1\kappa(u) < \epsilon_2\kappa(v)$.
    \end{enumerate}
\end{definition}

Essentially, for $\kappa\in\mathcal {K}(\alpha,\omega) $ if $\alpha(\epsilon_1 u) > \alpha(\epsilon_2 v)$ then $\epsilon_1\kappa( u) \geq \epsilon_2\kappa( v)$ and this inequality is strict if $\omega$ agrees with $\alpha$ here. See Example \ref{linear_extension_example}.

It will be useful later to know that if $(u_i)$ is a labeling of the vertices such that $|\alpha(u_1)| < |\alpha(u_2)| < \dots < |\alpha(u_n)|$, then $\mathcal{K}(\alpha,\omega)$ is the set of all colorings such that 
\begin{enumerate}
    \item If $i<j$, then $|\kappa(u_i)| \leq |\kappa(u_j)|$.
    \item $\sgn_\kappa(u_i) = \sgn_{\alpha}(u_i)$ or $\kappa(u_i)=0$ for all $i$.
    \item If $i<j$ and $\omega(\sgn_{\alpha}(u_i)  u_i) < \omega(\sgn_{\alpha}(u_j) u_j)$, then $\sgn_{\alpha}(u_i)\kappa(u_i) < \sgn_{\alpha}(u_j)\kappa(u_j)$.
\end{enumerate}

\noindent We can easily verify that for any $\omega\in \mathcal{L}(P)$,
\[
   \mathcal{C}(P) = \bigsqcup_{\alpha\in \mathcal{L}(P)}\mathcal{K}(\alpha,\omega).
\]
First we show $\mathcal K(\alpha_1,\omega)$ and $\mathcal K(\alpha_2,\omega)$ must be disjoint when $\alpha_1\neq \alpha_2$.
When $\alpha_1\neq\alpha_2$, there must be some $\epsilon_1 u$ and $\epsilon_2 v$ in $\overline{\Sigma}$ such that $\alpha_1(\epsilon_1 u) > \alpha_1(\epsilon_2 v)$ but $\alpha_2(\epsilon_1 u) < \alpha_2(\epsilon_2v)$.
Without loss of generality let $\omega(\epsilon_1u) > \omega(\epsilon_2 v)$.
Then any $\kappa\in\mathcal{K}(\alpha_1,\omega)\cap\mathcal{K}(\alpha_2,\omega)$ must satisfy both $\kappa(\epsilon_1u) > \kappa(\epsilon_2v)$  by (4.6.3) and $\kappa(\epsilon_1 u) \leq \kappa(\epsilon_2 v)$ by (4.6.1) and (4.6.2).

Now we show that for any $\omega\in\mathcal{L}(P)$ and for any $\kappa$ which preserves $P$, there exists $\alpha\in\mathcal{L}(P)$ such that $\kappa$ is contained in some $\mathcal{K}(\alpha,\omega)$.
Simply construct $\alpha$ such that for all $u$ and $v$ in $\Sigma$, if $|\kappa(u)| < |\kappa(v)|$ then $|\alpha(+u)| < |\alpha(+v)|$ and if $|\kappa(u)| = |\kappa(v)|$ then $\alpha$ does the opposite of what $\omega$ does i.e. $|\alpha(+u)| > |\alpha(+v)|$ iff $|\omega(+u)|<|\omega(+v)|$.
Also put $\sgn_\alpha(v) = \sgn_\kappa(v)$ when $\kappa(v)\neq 0$ and if $\kappa(v) = 0$ then put $\sgn_\alpha(v) = -\sgn_{\omega}(v)$.

\begin{example}\label{linear_extension_example}
Consider the following signed poset $P$ and its (oriented) covering graph
\[
\begin{tikzpicture}
    \node[vertex, label = above:$a$] (v1) at (0,0) {};
    \node[vertex, label = below:$b$] (v2) at (2,0) {};
    \node[vertex, label = above:$c$] (v3) at (1,1.712) {};

    \ltr {v1} {v2} [] [below];
    \ntp {v1} {v3} [] [above] [sloped];
    \ptn {v2} {v3} [] [above] [sloped];
    \draw[latex-latex, line width=0.5mm, loop below, out = 175, in = 250, looseness = 10] (v1) to (v1);
    \node[rotate = -45] (sign) at (-0.85,-0.6) {$-$};

     \node[vertex, label = above:$+c$] (+c) at (6.2,2) {};
        \node[vertex, label = above:$+b$] (+b) at (10.2,2) {};
        \node[vertex, label = left:$+a$] (+a) at (8.2,1) {};
        \node[vertex, label = below:$-c$] (-c) at (9.2,-1) {};
        \node[vertex, label = below:$-a$] (-a) at (7.2,0) {};
        \node[vertex, label = below:$-b$] (-b) at (5.2,-1) {};
        \drct {-b} {+c} {-latex};
        \drct {-a} {+c} {latex-};
        \drct {-a} {+a} {-latex};
        \drct {-c} {+a} {latex-};
        \drct {+a} {+b} {-latex};
        \drct {-b} {-a} {-latex};
        \drct {-c} {+b} {-latex};
\end{tikzpicture}
\]
this oriented covering graph has four possible linear extensions, which we will call $\alpha,\beta,\gamma$ and $\omega$:
\[
\begin{tikzpicture}

\node[label = right:$\alpha$] (a) at (0,4.5) {};
        \node[vertex, label = right:$+b$] (11) at (0,3.5) {};
        \node[vertex, label = right:$+a$] (12) at (0,2.5) {};
        \node[vertex, label = right:$-c$] (13) at (0,1.5) {};
        \node[vertex, label = right:$+c$] (14) at (0,0.5) {};
        \node[vertex, label = right:$-a$] (15) at (0,-0.5) {};
        \node[vertex, label = right:$-b$] (16) at (0,-1.5) {};
        \drct {12} {11} {-latex};
        \drct {13} {12} {-latex};
        \drct {14} {13} {-latex};
        \drct {15} {14} {-latex};
        \drct {16} {15} {-latex};

\node[label = right:$\beta$] (b) at (4,4.5) {};
        \node[vertex, label = right:$+c$] (21) at (4,3.5) {};
        \node[vertex, label = right:$+b$] (22) at (4,2.5) {};
        \node[vertex, label = right:$+a$] (23) at (4,1.5) {};
        \node[vertex, label = right:$-a$] (24) at (4,0.5) {};
        \node[vertex, label = right:$-b$] (25) at (4,-0.5) {};
        \node[vertex, label = right:$-c$] (26) at (4,-1.5) {};
        \drct {22} {21} {-latex};
        \drct {23} {22} {-latex};
        \drct {24} {23} {-latex};
        \drct {25} {24} {-latex};
        \drct {26} {25} {-latex};

\node[label = right:$\gamma$] (c) at (8,4.5) {};
        \node[vertex, label = right:$+b$] (31) at (8,3.5) {};
        \node[vertex, label = right:$+c$] (32) at (8,2.5) {};
        \node[vertex, label = right:$+a$] (33) at (8,1.5) {};
        \node[vertex, label = right:$-a$] (34) at (8,0.5) {};
        \node[vertex, label = right:$-c$] (35) at (8,-0.5) {};
        \node[vertex, label = right:$-b$] (36) at (8,-1.5) {};
        \drct {32} {31} {-latex};
        \drct {33} {32} {-latex};
        \drct {34} {33} {-latex};
        \drct {35} {34} {-latex};
        \drct {36} {35} {-latex};

\node[label = right:$\omega$] (w) at (12,4.5) {};
        \node[vertex, label = right:$+b$] (41) at (12,3.5) {};
        \node[vertex, label = right:$+a$] (42) at (12,2.5) {};
        \node[vertex, label = right:$+c$] (43) at (12,1.5) {};
        \node[vertex, label = right:$-c$] (44) at (12,0.5) {};
        \node[vertex, label = right:$-a$] (45) at (12,-0.5) {};
        \node[vertex, label = right:$-b$] (46) at (12,-1.5) {};
        \drct {42} {41} {-latex};
        \drct {43} {42} {-latex};
        \drct {44} {43} {-latex};
        \drct {45} {44} {-latex};
        \drct {46} {45} {-latex};
\end{tikzpicture}
\]
Then we can calculate 
\begin{itemize}
   %  \item $\mathcal{K}(\alpha,\omega)$ is the set of all colorings $\kappa$ with $|\kappa(c)| \geq |\kappa(b)| \geq |\kappa(a)|$, and $\kappa(c),\kappa(b),\kappa(a)\geq 0$ and also the following relations 
   %  \begin{align*}
   %      -&\kappa(b) <-\kappa(c), & &-\kappa(b)<-\kappa(a), &  &-\kappa(b)<\kappa(a), & &-\kappa(b) < \kappa(c) , & &-\kappa(b) < \kappa(b),  \\
   %      -&\kappa(c)<\kappa(a), &  &-\kappa(c)<\kappa(c), &  &-\kappa(c)<\kappa(b), &  & & &\\
   %      -&\kappa(a)<\kappa(a), & &-\kappa(a) < \kappa(c), &  & -\kappa(a)<\kappa(b), & & & &\\
   %      &\kappa(a)<\kappa(b), & & & & & &\\
   %      &\kappa(c)<\kappa(b) & & & & & & 
   %  \end{align*}
   % It turns out that most of these relations are redundant, and the only relation which gives us more information is $\kappa(c)<\kappa(b)$. 
   
   % This means that $\mathcal{K}(\alpha,\omega) = \{\kappa:\kappa(b)>\kappa(c)\geq \kappa(a)>0\}$.
   \item $\mathcal{K}(\alpha,\omega)$ is the set of all colorings $\kappa$ with $|\kappa(b)| \geq |\kappa(a)| \geq |\kappa(c)|$, $\kappa(b),\kappa(a)\geq 0$, $\kappa(c)\leq 0$ and also the following relations 
    \begin{align*}
        -&\kappa(b) <-\kappa(a),& -&\kappa(b)<\kappa(c), &  &-\kappa(b)<-\kappa(c), & &-\kappa(b) < \kappa(a), & &-\kappa(b) < \kappa(b),  \\
        -&\kappa(a) < \kappa(c),&   -&\kappa(a)<-\kappa(c),&  &-\kappa(a)<\kappa(a), & & -\kappa(a)<\kappa(b),  & & \\
        &\kappa(c)<\kappa(a),&  &\kappa(c) < \kappa(b),  &  &                     & &                        & &\\
        -&\kappa(c)<\kappa(a),& -&\kappa(c) < \kappa(b)  &  &                    & &                        & &\\
        &\kappa(a)<\kappa(b) & &                         &  &                   & &                        & & 
    \end{align*}
   It turns out that most of these relations are redundant, and the only relations which give us more information are $\kappa(a)<\kappa(b)$ and $-\kappa(c)<\kappa(a)$. 
   
   This means that $\mathcal{K}(\alpha,\omega) = \{\kappa:\kappa(b)>\kappa(a)> -\kappa(c)\geq 0\}$.
   \item We can note that the colorings of $\mathcal{K}(\beta, \omega)$ are those with $|\kappa(c)| \geq |\kappa(b)| \geq |\kappa(a)|$, $\kappa(b),\kappa(b),\kappa(a)\geq 0$, and the extra relations (since all of the others are be redundant) that $\kappa(a)<\kappa(b)$ and $-\kappa(a) < \kappa(a)$. So $\mathcal{K}(\beta,\omega) = \{\kappa: \kappa(c) \geq \kappa(b) > \kappa(a) > 0\} $.
   \item We may similarly calculate that $\mathcal{K}(\gamma,\omega) = \{\kappa:\kappa(b)>\kappa(c)\geq \kappa(a)>0\} $.
   \item Lastly, $\mathcal{K}(\omega,\omega) =  \{\kappa: \kappa(b) > \kappa(a) > \kappa(c) > 0\} $
\end{itemize}
We can note that these sets are pairwise disjoint, each coloring contained in one of these sets preserves $P$, and every coloring which preserves $P$ is contained in one of these sets.
\end{example}

For $\alpha,\omega\in \mathcal{L}(P)$, let us put $F_{\alpha,\omega} = \sum_{\kappa\in \mathcal{K}(\alpha,\omega)}x^{\kappa}$ so that for any fixed $\omega\in \mathcal{L}(P)$ we have $Y_P = \sum_{\alpha\in \mathcal{L}(P)}F_{\alpha,\omega}$. Then for any fixed $\omega\in \mathcal{L}(P)$ we have that
$$
    X_{\Sigma} = \sum_{P\in \mathcal{A}(\Sigma)}Y_P =  \sum_{P\in \mathcal{A}(\Sigma)}\left(\sum_{\alpha\in \mathcal{L}(P)}F_{\alpha,\omega}\right).
$$

%% file: Sections/sink_counting_appendix.tex
\section{A sink counting function}
% \begin{enumerate}
%     \item $F_{\alpha,\omega}$ is linearly independent for fixed $\omega$
%     \item Definition of $\phi_{\omega}$, it is well defined and works the way we want
%     \item construction of general $\phi$
%     \item $\phi$ is multiplicative
% \end{enumerate}

We are considering $F_{\alpha,\omega}$'s so that we may define a linear function, $\phi$, with the property that $\phi(Y_P) = t^{\sink(P)}$. Here $\sink(P)$ denotes the number of vertices which are sinks under the orientation $P$.
\subsection{Properties of $F_{\alpha,\omega}$}
First we need some preliminary facts about the $F_{\alpha,\omega}$'s. Let $\varnothing_n$ denote the trivial orientation on the signed graph with $n$ vertices and no edges, so that $\mathcal{L}(\varnothing_n)$ contains the linear extensions of any signed graph with $n$ vertices. For fixed $\omega\in \mathcal{L}(\varnothing_n)$, the set $\set{F_{\alpha,\omega}\st \alpha\in \mathcal{L}(\varnothing_n)}$ is linearly independent. Before this fact is proven, it should be noted that for fixed $\omega$, there may exist $\alpha\neq \beta$ such that $F_{\alpha,\omega} = F_{\beta,\omega}$. 

\begin{example}
    Let $\omega\in \mathcal{L}(\varnothing_2)$ such that $\omega(+v_1) = 1$ and $\omega(+v_2) = 2$, take $\alpha$ such that $\alpha(+v_1) = -1$ and $\alpha(+v_2) = 2$, and take $\beta$ such that $\beta(+v_2) = -1$ and $\beta(+v_1) = 1$. Then it can be computed that $F_{\alpha,\omega} = \sum_{0\leq i_1<i_2}x_{-i_1}x_{i_2} = F_{\beta,\omega}$.
\end{example}

We will also require an algebraic fact about the $F_{\alpha,\omega}$'s. We say that a power series with non-negative coefficients $A$ is a \emph{sub-sum} of a power series with non-negative coefficients $B$ if $B-A$ has non-negative coefficients. Equivalently, this means that every term which appears in $A$ also appears in $B$, including multiplicity. It is easily seen that if $A$ is a sub-sum of $B$ and $B$ is a sub-sum of $C$, then $A$ is a sub-sum of $C$ and additionally, if $A$ is a sub-sum of $B$ and $B$ is a sub-sum of $A$, then $A=B$.
\begin{lemma}
    If $F_{\alpha,\omega}$ is a sub-sum of $\sum_{i=1}^kF_{\beta_i,\omega}$, then $F_{\alpha,\omega}$ is a sub-sum of $F_{\beta_i,\omega}$ for some $i$.
\end{lemma} 
\begin{proof} 
    To see this consider a term $x_{a_1}^{b_1}x_{a_2}^{b_2}\dots x_{a_k}^{b_k}$ in $F_{\alpha,\omega}$ which has minimal length $k$ and with $0<a_1<a_2<\dots < a_k$. This term corresponds to a coloring in $\mathcal{K}(\alpha,\omega)$, call it $\kappa$, for which the least number of colors are used. This term must also exist on the right hand side, and so there is some $F_{\beta_i,\omega}$ which contains $x_{a_1}^{b_1}x_{a_2}^{b_2}\dots x_{a_k}^{b_k}$, coming from a coloring $\kappa'\in \mathcal{K}(\beta_i,\omega)$.

    Let $(v_j)$ be a labeling of the vertices such that $|\alpha(v_1)|<\dots< |\alpha(v_n)|$ and let $(u_j)$ be a labeling of the vertices such that $|\beta_i(u_1)|<\dots< |\beta_i(u_n)|$. Let $\theta$ be the function such that $\theta(v_j) = u_j$ for all $j$. Then $\kappa' = \kappa\circ \theta$ and moreover, by considering condition (3) in the definition of $\mathcal{K}(\alpha,\omega)$, we can see that any other coloring in $\mathcal{K}(\alpha,\omega)$ will become a coloring in $\mathcal{K}(\beta_i,\omega)$ when pre-composed with $\theta$. This is because any other coloring in $\mathcal{K}(\alpha,\omega)$ will use at least as many colors as $\kappa$ does and so when another coloring is pre-composed with $\theta$, it will satisfy the conclusion of condition (3) in the definition of $\mathcal{K}(\beta_i,\omega)$ whenever $\kappa\circ \theta$ does. Therefore $F_{\alpha,\omega}$ is a sub-sum of $F_{\beta_i,\omega}$.
\end{proof}

\begin{lemma}
For fixed $\omega\in \mathcal{L}(\varnothing_n)$, the set $\{F_{\alpha,\omega}\st \alpha\in \mathcal{L}(\varnothing_n)\}$ is linearly independent over $\mathbb{Q}$.
\end{lemma}
\begin{proof}
    Suppose that for some $\omega$, this set is linearly dependent. Then some non-trivial linear combination of elements of this set is equal to zero. After multiplying by constants, rearranging terms and writing terms of the form $n\cdot F_{\alpha,\omega}$ as $\sum_{k=1}^nF_{\alpha,\omega}$, the linear dependence equation becomes $\sum_{i=1}^{m}F_{\gamma_i,\omega} = \sum_{i=1}^{l}F_{\beta_i,\omega}$. We may also suppose that this linear dependence is minimal, i.e. all equal terms have been canceled and so $F_{\beta_i,\omega}\neq F_{\gamma_j,\omega}$ for all $i,j$.
    
    Pick a term on the right hand side which is not a sub-sum of any other term on the right hand side, except for terms it is equal to. There is a term like this because there are only finitely many terms. Without loss of generality, $F_{\beta_1,\omega}$ is such a term. $F_{\beta_1,\omega}$ is a sub-sum of the left hand side and so it is a sub-sum of a particular term. Without loss of generality, $F_{\beta_1,\omega}$ is a sub-sum of $F_{\gamma_1,\omega}$. $F_{\gamma_1,\omega}$ is a sub-sum of the right hand side and so it is a sub-sum of a particular term. Without loss of generality, $F_{\gamma_1,\omega}$ is a sub-sum of $F_{\beta_2,\omega}$. But this means that $F_{\beta_1,\omega}$ is also a sub-sum of $F_{\beta_2,\omega}$ and therefore $F_{\beta_1,\omega}=F_{\beta_2,\omega}$ by assumption. Then it must be that $F_{\beta_1,\omega} = F_{\gamma_1,\omega}$ since each is a sub-sum of the other. This is a contradiction since we assumed that this linear dependence was minimal. Therefore $\{F_{\alpha,\omega}\st \alpha\in \mathcal{L}(\varnothing_n)\}$ is linearly dependent for each $\omega\in \mathcal{L}(\varnothing_n)$.
\end{proof}
\subsection{Defining auxiliary functions}
Define $V_{\omega} = \spn\set{F_{\alpha,\omega}\st \alpha\in\mathcal{L}(\varnothing_n)}$.
For fixed $\omega$, we will define a linear function $\phi_{\omega}:V_{\omega} \rightarrow \Q[t]$. Given $\alpha\in\mathcal{L}(\varnothing_n)$, let $(v_i)$ be a relabeling of the vertices such that $|\alpha(+v_1)|<\dots <|\alpha(+v_n)|$. Also, let $\epsilon_i = \sgn_\alpha(v_i)$. Then define $\phi_{\omega}$ such that
\[
    \phi_{\omega}(F_{\alpha,\omega}) = \begin{cases}
        t(t - 1)^k &\text{ if } \omega(\epsilon_iv_i)>\omega(\epsilon_jv_j) \text{ for all } i,j \text{ with } n-k \leq i<j \leq n,\\
        & 0<\omega(\epsilon_iv_i)<\omega(\epsilon_jv_j) \text{ for all } i,j \text{ with } 1 \leq i<j \leq n-k,\\ 
        &\text{ and }\epsilon_i = + \text{ for } n-k\leq i \leq n, \text{ for } 0\leq k < n.\\
        (t - 1)^k &\text{ if } \omega(\epsilon_iv_i)>\omega(\epsilon_jv_j) \text{ for all } i,j \text{ with } n-k \leq i<j \leq n,\\
        & 0<\omega(\epsilon_iv_i)<\omega(\epsilon_jv_j) \text{ for all } i,j \text{ with } 1 \leq i<j \leq n-k,\\ 
        &\epsilon_i = + \text{ for } n-k< i \leq n \text{ and } \epsilon_{n-k} = -, \text{ for } 0\leq k < n.\\
        (t - 1)^n &\text{if } 0 > \omega(+v_{1}) > \omega(+v_{2}) > \dots > \omega(+v_{n})  \text{ and } \epsilon_i = + \text{ for all } i \\
        0 &\text{otherwise.}
    \end{cases}
\]
    It should be noted that the first case corresponds to the situation where largest element under $\omega$ is positive, namely it is $+v_{n-k}$, where $\alpha$ places $k-1$ positive elements above $+v_{n-k}$ and where the ordering from $\epsilon_1v_1$ to $+v_{n-k}$ agrees with $\omega$, i.e. $\alpha(\epsilon_1 v_1)<\alpha(\epsilon_2v_2)$ and $\omega(\epsilon_1 v_1)<\alpha(\epsilon_2v_2)$ and so on.
    
    The second case corresponds to an almost identical situation to the first case, except that the maximal element under $\omega$ is negative. The third case is the natural interpretation of the second case (i.e. $v_0=0$) with $k=n$.

    Currently it is not clear that $\phi_{\omega}$ is well defined. To show that it is, it suffices to show that whenever we have $F_{\alpha,\omega} = F_{\beta,\omega}$, we also have $\phi_{\omega}(F_{\alpha,\omega}) = \phi_{\omega}(F_{\beta,\omega})$.
    
    Suppose we have $F_{\alpha,\omega} = F_{\beta,\omega}$ for some $\alpha,\beta\in \mathcal{L}(\varnothing_n)$. Let $(v_i)$ be a labeling of the vertices such that $|\alpha(+ v_1)| < \dots <|\alpha(+v_n)|$ and let $\epsilon_i = \sgn_{\alpha}(v_i)$ for each $i$. 
    
    Let $(u_i)$ be a analogous relabeling of the vertices such that $|\beta(+u_1)|<\dots<|\beta(+u_n)|$. Since $F_{\alpha,\omega} = F_{\beta,\omega}$, we have that $\sgn_{\beta}(u_i) = \sgn_{\alpha}(v_i) = \epsilon_i$ for all $i$. 
   
   It must also be that any coloring in $\mathcal{K}(\alpha,\omega)$ when pre-composed with the map $u_j \mapsto v_j$ becomes a coloring in $\mathcal{K}(\alpha,\omega)$ and moreover, this is serves as a bijection between $\mathcal{K}(\alpha,\omega)$ and $\mathcal{K}(\beta, \omega)$. 
   
    Pick any $i$ and $j$ with $i<j$. We want to show that $\omega(\epsilon_iv_i)>\omega(\epsilon_jv_j)$ iff $\omega(\epsilon_iu_i)>\omega(\epsilon_ju_j)$.
   
    Suppose that $\omega(\epsilon_iv_i)>\omega(\epsilon_jv_j)$, then there is a $\kappa\in \mathcal{K}(\alpha,\omega)$ with $\epsilon_i\kappa(v_i) = \epsilon_j\kappa(v_j)$. This also means that there is a $\kappa'\in \mathcal{K}(\beta,\omega)$ with $\epsilon_i\kappa'(v_i) = \epsilon_j\kappa'(v_j)$ and hence that $\omega(\epsilon_iu_i)>\omega(\epsilon_ju_j)$. By symmetry we have the other direction, so that $\omega(\epsilon_iv_i)>\omega(\epsilon_jv_j)$ iff $\omega(\epsilon_iu_i)>\omega(\epsilon_ju_j)$ and therefore $\omega(\epsilon_iv_i)<\omega(\epsilon_jv_j)$ iff $\omega(\epsilon_iu_i)<\omega(\epsilon_ju_j)$. This fact makes it clear that $\phi_{\omega}(F_{\alpha,\omega}) = \phi_{\omega}(F_{\beta,\omega})$.

\subsection{Sink counting}
\begin{lemma}\label{lem:sink-counting}
    If $P$ is a signed poset and $\omega$ any linear extension of $P$, then $\phi_{\omega}(Y_P) = t^{\sink(P)}$.
\end{lemma}
\begin{proof}
    % Given some poset $P$, fix $\omega$ to be a linear extension of $P$.
    % 
    This proof is largely based on the proof of Theorem 3.3 in \cite{stanley95}.

    We will prove the lemma in each of the cases in the definition of $\phi$. We may do this since if $\alpha$ is a linear extension of $P$ such that $F_{\alpha,\omega}$ falls into the first case of $\phi_{\omega}$, then $F_{\beta,\omega}$ will fall into the first case or the zero case of $\phi_{\omega}$ for all other $\beta$ which are linear extensions of $P$.
    % In the first case of $\phi$, the maximal element under $\omega$ is positive, in the second case the maximal element under $\omega$ is negative with $k<n$, and the third case is equivalent to the second case if $k=n$ and $v_0$ is taken to be $0$. It is important to note that the vertex $v_j$ in $P$ is a sink iff $+v_j$ is the maximal element of some linear extension of $P$.
    
    First, consider the case where the largest element under $\omega$ is positive, call it $+s$, and so the vertex $s$ is a sink in $P$. Now select any $k$ of the remaining sinks of $P$ other than $s$ (of which there are $\sink(P)-1$ to choose from), and call these vertices $u_1, u_2,\dots, u_k$.
    
    The remaining $n-k-1$ vertices will be labeled $v_1,v_2,\dots, v_{n-k-1}$ such that $\abs{\omega(v_i)} < \abs{\omega(v_j)}$ iff $i<j$.
    
    Now, consider the linear extension $\alpha$, which puts
    \begin{enumerate}
        \item $\alpha(u_i) = n-i+1$,
        \item $\alpha(s) = n-k$,
        \item for $i<n-k$, $\abs{\alpha(v_i)} = i$ and $\sgn_\alpha(v_i) = \sgn_\omega(v_1)$.
    \end{enumerate}
    Hence this definition places $+u_1$ as the largest element under $\alpha$, $+u_2$ as the second largest, and so on. Then $+s$ is the greatest element below $+u_k$, and after $+s$, the vertices are arranged as their ordering in $\omega$ and each with the same sign as in $\omega$.
    
    Note that after initially choosing $k$ sinks no more choices are made, meaning that any $\alpha$ constructed this way is unique when given a choice of $k$ sinks.
    
    To see that $\alpha$ is a linear extension of $P$, we will examine every possible pair of vertices and determine that $\alpha$ respects the relation between them, if present in $P$.
    
    For any two vertices $v_i,v_j\in\set{v_1,\dots,v_{n-k-1}}$, it follows that $\abs{\alpha(v_i)} < \abs{\alpha(v_j)}$ iff $\abs{\omega(v_i)} < \abs{\omega(v_j)}$, with the signs of these vertices under $\alpha$ being the same as under $\omega$. Since $\omega$ respects all relations of $P$, this means $\alpha$ respects any relations between $v_i$ and $v_j$ present in $P$.
    
    So too $\alpha$ respects relations between any $p,q\in\set{u_1,\dots,u_k}\cup\set{s}$. Since they are all sinks of $P$, there is no directed positive edge or inward facing directed negative edge between any two of them. If there is an outward facing negative edge between two of them, this translates to the relation $p >_P -q$. But this relation also holds in $\alpha$ since $\alpha(+p) > 0$ for any $p\in\set{u_1,\dots, u_k}\cup\set{s}$. This logic holds for any two vertices that were sinks in $P$, not just elements $u_1,\dots,u_k,s$. So, the only consideration left is that of a sink and non-sink.
    
    Call the sink $p$ and the non-sink $r$. Then the edges possibly present in $P$ are a directed positive edge from $r$ to $p$ or an inward facing directed negative edge. These directed edges invoke the relations $p >_P r$ and $p >_P -r$, respectively. In either case we can see that $\alpha$ satisfies these relations as well since $\abs{\alpha(+p)} > \abs{\alpha(+r)}$ by construction.
    
    Therefore, $\alpha$ is a linear extension of $P$.
    
    So, for any $k<\sink(P)$ and for each choice of $k$ sinks of $P$ (other than $s$), there is exactly one linear extension $\alpha$ which satisfies the first case of $\phi$. This means that there are $N=\ncr{\sink{P}-1}{k}$ linear extensions $\alpha_1,\dots,\alpha_N$ for which $\phi(F_{\alpha_i,\omega}) = t(t-1)^k$. This holds for all $k\leq\sink{P}-1$ and any $\alpha$ not of this form must have $\phi(F_{\alpha,\omega}) = 0$ by uniqueness of $\alpha$'s.
    
    Hence we obtain
    \begin{align*}
        \phi(Y_P) &= \phi\left(\sum_{\alpha\in \mathcal{L}(P)}F_{\alpha,\omega}\right) = \sum_{k=0}^{\sink{P}-1}\ncr{\sink{P}-1}{k}t(t-1)^k \\
        &= t\cdot\sum_{k=0}^{\sink{P}-1}\ncr{\sink{P}-1}{k}(t-1)^k 
        = t\cdot t^{\sink{P}-1} = t^{\sink{P}}.
    \end{align*}
    %The second to last equality easily follows from the binomial theorem.
    This proves the theorem in the case where the maximal element under $\omega$ is positive.

    In the case where the maximal element under $\omega$ is negative the argument is almost identical, except that there are $\sink{P}$ vertices to choose from when selecting $u_1,\dots, u_k$, since $s$ is a source in $P$. The verification that the previous construction still works is straightforward. The third case of $\phi$ occurs when $\sink{P}=n$ and can be included in the construction.
    
    So, for any $k\leq\sink{P}$, and for each choice of $k$ sinks of $P$, there is exactly one linear extension $\alpha$ which satisfies the second or third case of $\phi$. Therefore there are $N=\ncr{\sink{P}}{k}$ linear extensions $\alpha_1,\dots,\alpha_N$ for which $\phi(F_{\alpha_i,\omega}) = (t-1)^k$. This holds for all $k\leq\sink{P}$ and any $\alpha$ not of this form must have $\phi(F_{\alpha,\omega}) = 0$ by uniqueness of $\alpha$'s.
    
    Hence we obtain
    \begin{align*}
        \phi(Y_P) &= \phi\left(\sum_{\alpha\in\mathcal{L}(P)}F_{\alpha,\omega}\right) \\
        &= \sum_{k=0}^{\sink{P}}\ncr{\sink{P}}{k}(t-1)^k 
        = t^{\sink{P}}.
    \end{align*}
    This completes the proof.
\end{proof}
\subsection{Constructing the sink counting function}

First we will define 
$$
\mathbb{Y} = \spn\set{Y_P\st P \text{ is a signed poset}}
$$
and note that $\mathbb{Y}$ is a subspace of $\bigoplus_{n=1}^{\infty} \sum_{\omega \in \mathcal{L} (\varnothing_n)} V_{\omega}$ since every $Y_P$ can be written as a sum of $F_{\alpha,\omega}$'s.

For fixed $n\in \N$ and for each $\omega\in \mathcal{L}(\varnothing_n)$, we have a linear function $\phi_{\omega}$ and moreover, these functions have the property that for any $\omega_1,\omega_2$, $\phi_{\omega_1}\mid_{V_{\omega_1}\cap V_{\omega_2}\cap \mathbb{Y}}=\phi_{\omega_2}\mid_{V_{\omega_1}\cap V_{\omega_2}\cap \mathbb{Y}}$. So we may define our desired function $\phi:\mathbb{Y}\rightarrow \Q[t]$ to agree with the $\phi_{\omega}$'s. More specifically, let $B$ be a basis for $\mathbb{Y}$ with the property that any element of $B$ is an element of $V_{\omega}$ for some $\omega$. Then for any $b\in B$, pick $\omega$ such that $b\in V_{\omega}$ and define $\phi(b) = \phi_{\omega}(b)$. Our choice of $\omega$ is not unique, but $\phi$ does not depend on choices of $\omega$ since $\phi_{\omega_1}$ and $\phi_{\omega_2}$ agree on $\mathbb{Y}$ when applicable. Therefore we have a well defined linear function $\phi:\mathbb{Y}\rightarrow \Q[t]$ such that $\phi(Y_P) = t^{\sink(P)}$.

%% file: Sections/phi_calcs.tex
\section{Using $\phi$ for calculations}
\begin{theorem}
    Let $\Sigma$ be a signed graph and let $\acyc_{\Sigma}(k)$ denote the number of acyclic orientations of $\Sigma$ which have $k$ sinks. Then $\phi(X_\Sigma)=\sum_{k=0}^\oo\acyc_\Sigma(k)t^k$.
\end{theorem}
\begin{proof}
    We know that $X_{\Sigma} =  \sum_{P\in\mathcal A(\Sigma)}Y_P$, so we may apply $\phi$ to arrive at
    \[
        \varphi(X_{\Sigma}) = \sum_{P\in\mathcal A(\Sigma)}\varphi(Y_P) = \sum_{P\in\mathcal A(\Sigma)}t^{\sink(P)}.
    \]
    After counting terms, we have $\sum_{P\in\mathcal A(\Sigma)}t^{\sink(P)} = \sum_{k=0}^{\infty}\acyc_{\Sigma}(k)t^k$ as desired.
\end{proof}

\begin{observation}
    If $P_1$ and $P_2$ are disjoint signed posets, then $Y_{P_1\sqcup P_2} = Y_{P_1}\cdot Y_{P_2}$. Since also $\sink(P_1\sqcup P_2) = \sink(P_1)+\sink(P_2)$, we can see that
    \[
        \phi(Y_{P_1})\cdot \phi(Y_{P_2}) = t^{\sink(P_1)}\cdot t^{\sink(P_2)}= t^{\sink(P_1)+\sink(P_2)} = \phi(Y_{P_1\sqcup P_2}) = \phi(Y_{P_1}\cdot Y_{P_2}).
    \]
    Since $\mathbb{Y}$ is spanned by the $Y_P$'s, this means that for arbitrary $f, g\in \mathbb{Y}$, $\phi(f\cdot g) = \phi(f)\cdot \phi(g)$.
\end{observation}

\begin{definition}
    For $k\geq 0$, let $S_k$ be the star graph which has $k$ positive edges and $k+1$ vertices.
\end{definition}

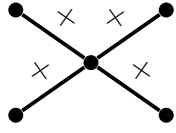
\begin{figure}[h!t]
    \centering
    \begin{tikzpicture}
        \node[small] (v1) at (0,0) {};
        \node[small] (v2) at (1,0.7) {};
        \node[small] (v3) at (1,-0.7) {};
        \node[small] (v4) at (-1,0.7) {};
        \node[small] (v5) at (-1,-0.7) {};
        \pstv {v1}{v2};
        \pstv {v1}{v3};
        \pstv {v1}{v4};
        \pstv {v1}{v5};
    \end{tikzpicture}
    \caption{The signed graph $S_4$.}
\end{figure}

\begin{lemma} $X_{S_k} = \sum_{i=0}^k(-1)^i\ncr ki p_{1,0}^{k-i}\,p_{i+1,0}$.
\end{lemma}
\begin{proof}
    Consider the terms obtained by contracting $i$ of the $k$ total edges and deleting the rest. This term will have a factor of $(-1)^i$ since the contracted term is always subtracted. When $i$ of the edges are contracted, there are $k+1-i$ vertices left. The one in the center has weight $(i+1,0)$, and the vertices around it have weight $(1,0)$. The term this corresponds to is $(-1)^ip_{1,0}^{k-i}p_{i+1,0}$. Finally, there are $\ncr ki$ ways to contract $i$ edges and we can do this for any $i=0,1,\dots, k$, demonstrating equality.
\end{proof}

\begin{theorem}\label{thm:ez-computation}
    For any $a\in\N$, $\phi(p_{a, 0}) = (t-1)^a-(-1)^a$.
\end{theorem}
\begin{proof}
    For any $k$, we have that $X_{S_k} = \sum_{i=0}^k(-1)^i\ncr ki p_{1,0}^{k-i}\,p_{i+1,0}$, where $S_k$ is the star graph with $k$ edges. To find $\phi(X_{S_k})$, consider the orientation of $S_k$ with $i$ edges pointing away from the center vertex. If $i\neq 0$, then there are $i$ sinks and if $i = 0$ then the vertex in the center is a sink. There are $\ncr ki$ orientations of this form and so
    \[
        \phi(X_{S_k}) = \sum_{i=0}^{k}\acyc_{\Sigma}(k)t^i = t+\sum_{i=1}^{k}\ncr kit^i = (t-1)+\sum_{i=0}^{k}\ncr kit^i = (t-1)+(t+1)^k.
    \]

    We proceed by induction. First, note that $\phi(p_{1,0}) = t$, since $p_{1,0}$ is the chromatic $B$-Symmetric polynomial of a single vertex with no edges. Suppose that $\phi(p_{a,0}) = (t-1)^a-(-1)^a$ for all positive $a\leq k$. We apply this assumption to $\phi(X_{S_k})$ to see
    \begin{align*}
        \phi(X_{S_k}) &= \sum_{i=0}^k(-1)^i\phi(p_{1,0}^{k-i})\cdot\phi(p_{i+1,0}) \\
        &= (-1)^k\phi(p_{k+1,0})+\sum_{i=0}^{k-1}(-1)^i\ncr kit^{k-i}\cdot\left((t-1)^{i+1}-(-1)^{i+1}\right) \\
        &= (-1)^k\phi(p_{k+1,0})+\sum_{i=0}^{k-1}(-1)^i\ncr kit^{k-i}(t-1)^{i+1}-\sum_{i=0}^{k-1}(-1)^i\ncr kit^{k-i}(-1)^{i+1} \\
        &= (-1)^k\phi(p_{k+1,0})+(t-1)\sum_{i=0}^{k-1}\ncr{k}{i}t^{k-i}(1-t)^{i}+\sum_{i=0}^{k-1}\ncr {k}{i}t^{k-i}\\
        &= (-1)^k\phi(p_{k+1,0})+ (t-1)\left(1-(-1)^{k}(t-1)^{k}\right)+\left((t+1)^k-1\right). 
    \end{align*}
    We know that $\phi(X_{S_k}) = (t-1)+(t+1)^k$, so we have that $(-1)^k\phi(p_{k+1,0})-(-1)^k(t-1)^{k+1}-1 = 0$ and hence $\phi(p_{k+1,0})=(t-1)^{k+1}-(-1)^{k+1}$. So, by induction we are done.
\end{proof}

\begin{lemma}
    Let $\S_{n,m}$ denote the signed graph created by connecting the center vertices of $S_n$ and $S_m$ with a negative edge, then
    \[
        X_{\S_{n,m}}=\sum_{i=0}^{n}\left(\sum_{j=0}^{m}\ncr{n}{i}\ncr{m}{j}(-1)^{i+j}p_{1,0}^{n+m-i-j}(p_{i+1,0}p_{j+1,0}-p_{i+1,j+1})\right).
    \]
\end{lemma}
\begin{figure}[h!t]
    \centering
    \begin{tikzpicture}
        \node[small] (a1) at (0,0) {};
        \node[small] (a2) at (1,1) {};
        \node[small] (a3) at (1,-1) {};
        \node[small] (c1) at (1,0) {};
        
        \node[small] (c2) at (3,0) {};
        \node[small] (b1) at (3.5,0.867) {};
        \node[small] (b2) at (3.5,-0.867) {};
        
        \pstv {c1} {a1} [] [below];
        \pstv {c1} {a2};
        \pstv {c1} {a3};
        \pstv {c2} {b1};
        \pstv {c2} {b2} [] [below];
        \ngtv {c1} {c2};
    \end{tikzpicture}
    \caption{The graph $\S_{2,3}=\S_{3,2}$.}
\end{figure}
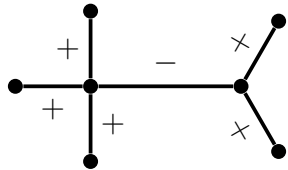
\begin{proof}
    We will perform weighted contraction deletion. Suppose that $i$ of the positive edges in the $S_n$ portion of the graph (of which there are $\ncr{n}{i}$ combinations) and $j$ of the positive edges in the $S_m$ portion of the graph (of which there are $\ncr{m}{j}$ combinations) are contracted, and the rest of the positive edges are deleted. This term will have a factor of $(-1)^{i+j}$ since contraction-terms are subtracted from deletion-terms. Furthermore, we are left with $(n-i)+(m-j)$ disjoint vertices of weight $(1,0)$ and two vertices of weights $(i+1,0)$ and $(j+1,0)$ which are connected by a negative edge. Performing weighted deletion-contraction on this last edge gives us the term
    \begin{align*}
        (-1)^{i+j}\left(p_{1,0}^{n+m-i-j}p_{i+1,0}p_{j+1,0}-p_{1,0}^{n+m-i-j}p_{i+1,j+1}\right)& \\
        &= (-1)^{i+j}p_{1,0}^{n+m-i-j}(p_{i+1,0}p_{j+1,0}-p_{i+1,j+1})
    \end{align*}
    There are $\ncr{n}{i}\ncr{m}{j}$ ways to arrive at this term, and so
    \[
        X_{\S_{n,m}} = \sum_{i=0}^{n}\left(\sum_{j=0}^{m}\ncr{n}{i}\ncr{m}{j}(-1)^{i+j}p_{1,0}^{n+m-i-j}(p_{i+1,0}p_{j+1,0}-p_{i+1,j+1})\right)
    \]
    as desired.
\end{proof}

\begin{theorem}\label{thm:ez-computation2}
    If $a,b\geq 1$, then $\phi(p_{a, b}) = -(-1)^{a+b}$.
\end{theorem}
\begin{proof}
    We will begin by noting that 
    \[
        X_{\S_{n,m}}-X_{S_n\sqcup S_m} = X_{\S_{n,m}}-X_{S_n}X_{S_m} = -\sum_{i=0}^n\sum_{j=0}^m\ncr{n}{i}\ncr{m}{j}(-1)^{i+j}p_{1,0}^{n+m-i-j}p_{i+1,j+1}
    \] We are able to calculate $\phi(X_{\S_{n,m}}-X_{S_n\sqcup S_m})$ by noting that the graph $S_n\sqcup S_m$ is almost identical to the graph $\S_{n,m}$, but $\S_{n,m}$ has a negative edge between the center vertices of the two star graphs. Any orientation of the graph $S_n\sqcup S_m$ will have the same number of sinks as the corresponding orientation of $\S_{n,m}$ for which the negative edge points outward. All orientations of both $S_n\sqcup S_m$ and $\S_{n,m}$ are acyclic and so it suffices to consider the orientations of $\S_{n,m}$ for which the negative edge points inward. Consider an orientation of $\S_{n,m}$ such the the negative edge points inward, with $i$ of the edges belonging to $S_n$ and $j$ of the edges belonging to $S_m$ pointing away from their center vertex. This orientation will have $i+j$ sinks and there are $\ncr{n}{i}\ncr{m}{j}$ orientations of this form. Therefore 
    \begin{align*}
        \phi(X_{\S_{n,m}}-X_{S_n\sqcup S_m})= \sum_{i=0}^n\sum_{j=0}^m\ncr{n}{i}\ncr{m}{j}t^{i+j} = (t+1)^{n+m}
    \end{align*}
    
    We will proceed by induction. We know that $\phi(p_{1,1}) = -1$ because $X_{\begin{tikzpicture}
        \node[small] (v) at (0,0) {};
        \node[small] (u) at (0.5,0) {};
        \edge {v} {u};
        \node (sign) at (0.25,0.125) {-};
    \end{tikzpicture}} = p_{1,0}^2-p_{1,1}$ and $\phi(X_{\begin{tikzpicture}
        \node[small] (v) at (0,0) {};
        \node[small] (u) at (0.5,0) {};
        \edge {v} {u};
        \node (sign) at (0.25,0.125) {-};
    \end{tikzpicture}}) = p_{1,0}^2-p_{1,1} = t^2+1$, so $\phi(p_{1,1})=-1$. Now suppose that $\phi(p_{a,b}) = (-1)^{a+b+1}$ for all $a,b\geq 1$ with $a+b\leq n+m+2$. Then by the induction hypothesis, we have that 
    \begin{align*}
        &\phi(X_{\S_{n,m}}-X_{S_n\sqcup S_m}) = \phi\left(-\sum_{i}^n\sum_{j}^m\ncr{n}{i}\ncr{m}{j}(-1)^{i+j}p_{1,0}^{n+m-i-j}p_{i+1,j+1}\right) \\
        &= -\phi(p_{n+1,m+1})(-1)^{n+m}-1-\sum_{i=0}^n\sum_{j=0}^m\ncr{n}{i}\ncr{m}{j}(-1)^{i+j}t^{n+m-i-j}(-1)^{(i+1)+(j+1)+1}\\
        &= -\phi(p_{n+1,m+1})(-1)^{n+m}-1+(t+1)^{n+m}
    \end{align*}
    We have that $\phi(X_{\S_{n,m}}-X_{S_n\sqcup S_m})= (t+1)^{n+m}$ and so it must be that $-\phi(p_{n+1,m+1})(-1)^{n+m}-1$ and so $\phi(p_{a+1,b+1}) = (-1)^{n+m+1} = (-1)^{(n+1)+(m+1)+1}$ which completes the proof.
\end{proof}

\begin{lemma}\label{thm:ez-computation3}
    $\phi(x_0) =-1$.
\end{lemma}
\begin{proof}
    Let $\Sigma$ be the signed graph 
    \begin{tikzpicture}
        \node[small] (v1) at (0,0) {};
    
        \draw[line width=0.5mm, loop above, out = 45, in = 135, looseness = 10] (v1) to (v1);
        \node (sign) at (0,0.5) {$-$};
    \end{tikzpicture}. Then $X_{\Sigma} = p_{1,0}-x_0$. $\Sigma$ has two orientations, one with 1 sink and one with no sinks and so $\phi(X_{\Sigma}) = t+1$. $\phi(p_{1,0}) =t$ and so it must be that $\phi(x_0) = -1$
\end{proof}
%     We will use the $B$-symmetric chromatic polynomial $X_{S_n}$ of the star graph $S_n$ with $n$ edges. Specifically, recall we have that $\phi(X_{S_n}) = (t+1)^n+(t-1)$. If we add a negative loop to the center edge of the star to obtain $S_n^-$, then it follows combinatorially that $\phi(X_{S_n^-}) = 2(t+1)^n+(t-1)$. The proof will proceed by induction.
    
%     First, $\phi(x_0^0) = 1 = (-1)^0$ easily. Now, assume that $\phi(x_0^n)=(-1)^n$ for all $0\leq n<k$ for some $k\in\N$, and consider the graph $S_k^-$. We can apply deletion-contraction to see that
%     \[
%         X_{S_k^-} = \sum_{i=0}^k(-1)^i\ncr ki p_{1,0}^{k-i}\,p_{i+1,0}-\sum_{i=0}^k(-1)^i\ncr ki p_{1,0}^{k-i}\left(-x_0^{i+1}\right).
%     \]
%     Notice that the first sum is precisely $Y_{S_k}$, and so we can apply $\phi$ and simplify the expression to
%     \begin{align*}
%         \phi(X_{S_k^-}) &= \phi(X_{S_k}) + \sum_{i=0}^{k-1}(-1)^i\ncr ki\phi\left(p_{1,0}^{k-i}\right)\phi\left(x_0^{i+1}\right)+(-1)^k\phi(x_0^k) \\
%         &= (t+1)^k+(t-1)-\sum_{i=0}^{k-1}(-1)^{2i+1}\ncr kit^{k-i}+(-1)^k\phi(x_0^k) \\
%         2(t+1)^k+(t-1) &= (t+1)^k+(t-1)+\sum_{i=0}^{k}\ncr kit^{k-i}-1+(-1)^k\phi(x_0^k) \\
%         1 &= (-1)^k\phi(x_0^k)
%     \end{align*}
    
%     Hence the theorem is proved.

Consider the \emph{elementary symmetric functions} in the variables $\dots,x_{-2}, x_{-1}, x_{0}, x_1,x_2, \dots$

$e_n = \sum_{i_1< i_2 < \dots < i_n}x_{i_1}x_{i_2}\cdots x_{i_n}$, where the sum ranges over all integer valued increasing sequences of length $n$.

\begin{lemma}
    $\phi(e_n) = t$ for all $n>0$.
\end{lemma}
\begin{proof}
    Newton's Identities state that $n\cdot e_n = (-1)^{n+1}p_{n,0}+\sum_{i=1}^{k-1}(-1)^{i+1}e_{k-i}p_{i,0}$.
    
    It is easy to verify that $\phi(e_1) = \phi(e_2) = t$, since $e_1 = p_{1,0}$ and $2e_2 = e_1p_{1,0}-p_{2,0}$. 
    
    We will proceed by induction. Suppose that for all $k\leq n-1$, we have that $\phi(e_k) = t$. Then 
    \begin{align*} 
        \phi(n\cdot e_n) &= (-1)^{n+1}\phi(p_{n,0})+\sum_{i=1}^{n-1}(-1)^{i+1}\phi(e_{k-i})\cdot\phi(p_{i,0})\\
        &= (-1)^{n+1}((t-1)^n-(-1)^n)+\sum_{i=1}^{n-1}(-1)^{i+1}t\cdot((t-1)^i-(-1)^i)\\
        &= 1-(1-t)^n+t\sum_{i=1}^{n-1}(1-(1-t)^i)\\
        &= 1-(1-t)^n+(n-1)t-t\sum_{i=1}^{n-1}(1-t)^i
    \end{align*}
    So in order to show that $\phi(e_n) = t$ for all $n$, it suffices to show that $1-(1-t)^n-t\sum_{i=1}^{n-1}(1-t)^i = t$ for all $n$, but this can easily be seen by induction since $1-(1-t)^2-t\sum_{i=1}^{2-1}(1-t)^i = 1-(1-2t+t^2)-t(1-t) = t$ and $1-(1-t)^n-t\sum_{i=1}^{n-1}(1-t)^i = 1-(1-t)^{n+1}-t\sum_{i=1}^{(n+1)-1}(1-t)^i$ for all $n$.
\end{proof}

%% file: Sections/main_theorem_proofs.tex
\begin{theorem}\label{main_thm_again}
    If the chromatic $B$-symmetric function $X_\Sigma$ of some signed graph $\Sigma$ is written in terms of sums and products from the augmented elementary $B$-symmetric basis, then the number of acyclic orientations of $\Sigma$ with $k$ sinks is the sum of the coefficients of terms having $k$ elementary $B$-symmetric function factors.
\end{theorem}

\begin{proof}
    We know that $\phi(q_{a,b})=\phi(z) = 1$ for all $a,b,n\geq 1$ and that $\phi(e_n) = t$ for all $n\geq 1$. We also know that the coefficient of $t^k$ in $\phi(X_{\Sigma})$ is the number of acyclic orientations of $\Sigma$ with $k$ sinks. The terms which $\phi$ sends to a multiple of $t^k$ are precisely those terms which have exactly $k$ elementary $B$-symmetric function factors.
\end{proof}

See Example \ref{main_thm_1_ex}.

\begin{definition}
    Let $\xi_n = \sum_{a=1}^n\ncr{n}{a}p_{a, 0}$ for $n\geq 1$. 
\end{definition}
Alternatively, we may write this as $p_{n,0} = \sum_{i=1}^n{\ncr{n}{i}}(-1)^{n-i}\xi_i$.
\begin{lemma}\label{main_lem_2}
    $\phi(\xi_n) = t^n$ for all $n$.
\end{lemma}
\begin{proof}
    \begin{align*}
        \phi(\xi_n) &= \sum_{a=1}^n\ncr{n}{a}\varphi(p_{a, 0}) \\
        &= \sum_{a=1}^n\ncr{n}{a}((t-1)^a-(-1)^a) \\
        &= \left(1 + \sum_{a=1}^n\ncr{n}{a}(t-1)^a\right)-\left(1+\sum_{a=1}^n\ncr{n}{a}(-1)^a\right) \\
        &= t^n-0.
    \end{align*}
\end{proof}

\begin{theorem}
    If the chromatic $B$-symmetric function $X_\Sigma$ of some signed graph $\Sigma$ is written in terms of sums and products from the set $\set{\xi_n\st n\geq 1}\cup\set{q_{a,b}\st a,b\geq 1}\cup\{z\}$, then the number of acyclic orientations of $\Sigma$ with $k$ sinks is the sum of the coefficients of terms such that the sum of the indices of each $\xi_n$ factor is equal to $k$.
\end{theorem}

This follows by Lemma \ref{main_lem_2} from the same argument as Theorem \ref{main_thm_again}. See Example \ref{main_thm_2_ex}.

%% file: Sections/further_results.tex
\section{further results}
As mentioned previously, a signed graph with all positive edges can be considered as an unsigned graph an vice versa. Suppose that $\Sigma$ is a signed graph with all positive edges. Let $|\Sigma|$ denote the unsigned graph which corresponds to $\Sigma$, i.e. $|\Sigma|$ would become $\Sigma$ if we added plus signs to each of its edges. We can note that the proper colorings of $|\Sigma|$ are precisely the proper colorings of $\Sigma$ which only use positive colors. So we may define $\proj_{>0 }$ to be a linear and multiplicative function such that $\proj_{>0}(x_i) = x_i$ for $i\geq 1$ and $\proj_{>0}(x_i) = 0$ for $i\leq 0$. It is clear that the chromatic $B$-symmetric function of $\Sigma$ becomes the chromatic symmetric function of $|\Sigma|$ when each $x_i$ which has a non-positive index is replaced by zero, i.e. $\proj_{>0}(X_{\Sigma}) = X_{|\Sigma|}$. It should be noted that due to the process of weighted deletion-contraction, $X_{\Sigma}$ will not contain any terms of the form $q_{a,b}$ or $z$ since $\Sigma$ has no negative edges. 

It is easy to see that $\proj_{>0}$ sends the elementary $B$-symmetric functions to the elementary symmetric functions. This means that Stanley's result about the elementary symmetric basis \cite[Theorem 3.3]{stanley95} follows immediately from the result about the augmented elementary $B$-symmetric basis (Theorem \ref{main_1}). 

Additionally, since $\proj_{>0}(\xi_n) = \sum_{a=1}^n{\ncr{n}{a}}p_{a}$, where $p_a = \sum_{i\geq 0}x_i^a$, we may define $\zeta_n = \sum_{a=1}^n{\ncr{n}{a}}p_{a}$ so that we have the following result. 
\begin{corollary}
If the chromatic symmetric function $X_G$ of some graph $G$ is written in terms of sums and products from the set $\set{\zeta_n\st n\geq 1}$, then the number of acyclic orientations of $G$ with $k$ sinks is the sum of the coefficients of terms such that the sum of the indices of each $\zeta_n$ factor is equal to $k$.
\end{corollary}
See Example \ref{main_thm_3_ex}.

We may also use $\phi$ to recover Zaslavsky's result \cite [Corollary 4.1]{zaslavsky}. Specifically, if $\chi_{\Sigma}$ is the signed chromatic polynomial of a signed graph $\Sigma$, i.e. $\chi_{\Sigma}(\lambda)$ is the number of proper colorings of $\Sigma$ in the colors $-\lambda,\dots, -1, 0, 1, \dots, \lambda$, then we may define $f:\mathbb{Y}\rightarrow \Q[\lambda]$ to be linear and multiplicative with $f(p_{a,b}) = 2\lambda+1$ for all $a,b\geq 0$ and $f(x_0) = 1$, Then $f$ has the property that $f(X_{\Sigma}) = \chi_{\Sigma}$. Next, note that $\phi(X_{\Sigma})\mid_{t=1}$ is equal to the total number of acyclic orientations of $\Sigma$ and that $\phi(p_{a,b})\mid_{t=1} = (-1)^{a+b+1}$ for all $a\geq 1$, $b\geq 0$. To recover Zaslavsky's result, we must show that $f(X_{\Sigma})\mid_{\lambda= -1} = (-1)^n\cdot(\text{The number of acyclic orientations of } \Sigma)$ when $\Sigma$ is a signed graph with $n$ vertices.

To do this, may note that $\phi(p_{a_1,b_1}\dots p_{a_k,b_k}x_0^c)\mid_{t=1} = (-1)^{n+k}$ where $n = a_1+b_1+\dots + a_k+b_k+c$ and $k$ is the number of $p_{a,b}$ terms. Also, we have that $f(p_{a_1,b_1}\dots p_{a_k,b_k}x_0^c)\mid_{\lambda=-1} = (-1)^k$ since $f(p_{a,b})\mid_{\lambda=-1} =-1$ and $f(x_0)\mid_{\lambda=-1} = 1$. Therefore if $\Sigma$ is a signed graph with $n$ vertices, then we have
$$
\chi_{\Sigma}(-1)=
f(X_{\Sigma})\mid_{\lambda= -1} = (-1)^n\phi(X_{\Sigma})\mid_{t=1} = (-1)^n\cdot(\text{\# of acyclic orientations of } \Sigma)
$$
which is Zaslavsky's result.

An equivalent way of stating this is that the sum of the absolute value of the coefficients of $X_\Sigma$ written in the $p$-basis is equal to the total number of acyclic orientations of $\Sigma$.

\begin{example}
    Let $\Sigma$ be the signed graph$\begin{tikzpicture}[baseline = 0.25cm]
        \node[small] (v1) at (0,0) {};
        \node[small] (v2) at (0.75,0) {};
        \node[small] (v3) at (0.375,0.65) {};
        \pstv {v1} {v2} [] [below];
        \ngtv {v1} {v3};
        \pstv {v2} {v3};
    \end{tikzpicture}$. 
    
    Then 
$$
        X_{\Sigma} = p_{1,0}^3-p_{1,0}p_{1,1}-2p_{1,0}p_{2,0}+2p_{2,1}+p_{3,0}-x_0^3
$$
    So $\Sigma$ has $|1|+|-1|+|-2|+|2|+|1|+|-1| = 8$ acyclic orientations.
\end{example}
We can also see that this holds for unsigned graphs as well

\begin{example}
    Let $G$ be the unsigned graph $\begin{tikzpicture}[baseline = 0.25cm]
        \node[small] (v1) at (0,0) {};
        \node[small] (v2) at (0.75,0) {};
        \node[small] (v3) at (0.375,0.65) {};
        \edge {v1} {v2} [] [below];
        \edge {v1} {v3};
        \edge {v2} {v3};
    \end{tikzpicture}$
    
    Then 
\[
    X_G = p_1^3-3p_1p_2+2p_3.
\]
 So, $G$ has $|1|+|-3|+|2|= 6$ acyclic orientations. 
\end{example}

%% file: Sections/bib.tex
\bibliographystyle{abbrv}

%% file: ms.bbl
\begin{thebibliography}{3}

    \bibitem{Chmutov_paper} S. Chmutov, S. Duzhin, S. Lando, \textit{Vassil knot invariants III. Forest algebra and weighted graphs}, Advances in Soviet Mathematics 21 (1994) 135-145.
    
    \bibitem{Chmutov_slides} S. Chmutov, \textit{B-Symmetric chromatic function of signed graphs}. \url{https://people.math.osu.edu/chmutov.1/talks/2020/slides-Moscow.pdf}.
    Video is available at \url{https://www.youtube.com/watch?v=khA7rP84sYY}.
    
    \bibitem{Crew} L. Crew, S. Spirkl, \textit{A deletion-contraction relation for the chromatic symmetric function}. European Journal of Combinatorics 89 (2020) 103-143.
    
    \bibitem{Egge} E. Egge, \textit{A Chromatic Symmetric Function for Signed Graphs}. \url{https://www.ericegge.net/slides/athens_slides.pdf}.
    
    \bibitem{Noble}  S. Noble, D. Welsh, \textit{A weighted graph polynomial from chromatic invariants of knots}, Annales de l’institut Fourier 49(3) (1999) 1057–1087.
    
    \bibitem{Rushil} R. Raghavan, \textit{A Symmetric Chromatic Function for Signed Graphs}.
    \url{https://people.math.osu.edu/chmutov.1/wor-gr-su19/Rushil-Raghavan-MIGHTY_10192019.pdf}
    
    \bibitem{stanley73} R. Stanley, \textit{Acyclic orientations of graphs}. Discrete Mathematics 5 (1973) 171-178.
    
    \bibitem{stanley95} R. Stanley, \textit{A symmetric function generalization of the chromatic polynomial of a graph}, Advances in Math. 111(1) (1995) 166–194.
    
    \bibitem{zaslavsky} T. Zaslavsky, \textit{Signed graph coloring}, Discrete Mathematics 39(2) (1982) 215–228.
\end{thebibliography}
